\documentclass[10 pt,a4paper]{article}

\usepackage{amsfonts,amsmath,amsthm,amssymb,graphicx}

\newtheorem{Proposizione}[subsection]{Proposition}
\newtheorem{Lemma}[subsection]{Lemma}
\newtheorem{Teorema}[subsection]{Theorem}
\newtheorem{Remark}[subsection]{Remark}
\newtheorem{Corollario}[subsection]{Corollary}

\numberwithin{equation}{section}
\begin{document}

\author{Bruno Bianchini \and Luciano Mari \and Marco Rigoli }
\title{\textbf{Spectral radius, index estimates\\
for Schr\"odinger operators\\
and geometric applications}}
\date{}
 \maketitle
\scriptsize \begin{center} Dipartimento di Matematica,
Universit\`a
degli Studi di Milano,\\
Via Saldini 50, I-20133 Milano (Italy)\\
E-mail addresses: bianchini@dmsa.unipd.it, lucio.mari@libero.it,
Marco.Rigoli@mat.unimi.it
\end{center}

\newcommand{\R}{\mathbb R}
\newcommand{\N}{\mathbb N}
\newcommand{\bull}{\rule{2.5mm}{2.5mm}\vskip 0.5 truecm}
\newcommand{\binomio}[2]{\genfrac{}{}{0pt}{}{#1}{#2}} 
\begin{abstract}
In this paper we study the existence of a first zero and the
oscillatory behavior of solutions of the ordinary differential
equation $(vz')'+Avz = 0$, where $A,v$ are functions arising from
geometry. In particular, we introduce a new technique to estimate
the distance between two consecutive zeros. These results are
applied in the setting of complete Riemannian manifolds: in
particular, we prove index bounds for certain Schr\"odinger
operators, and an estimate of the growth of the spectral radius of
the Laplacian outside compact sets when the volume growth is
faster than exponential. Applications to the geometry of complete
minimal hypersurfaces of Euclidean space, to minimal surfaces and
to the Yamabe problem are discussed.
\end{abstract}
\normalsize
\section{Introduction}
Radialization techniques are a powerful tool in investigating
complete Riemannian manifolds. In favourable circumstances these
lead to the study of an ordinary differential equation in order to
control the solutions of a given partial differential equation. In
this respect, one of the challenging problems involved is the
study of the sign of the solutions of the ODE, and the positioning
of the possible zeros. In this paper we determine some conditions
ensuring the oscillatory behavior, the existence of zeros and
their positioning, of a solution $z(t)$ of the following Cauchy
problem:
\begin{equation}\label{cauchy}
\left\{ \begin{array}{l} (v(t)z'(t))' + A(t)v(t) z(t) = 0 \qquad
\mathrm{on \ }
(0,+\infty)\\[0.1cm]
z'(t) = O(1) \quad \mathrm{as \ } t \downarrow 0^+ \quad ,
\quad z(0^+) = z_0>0
\end{array}\right.
\end{equation}
where $v(t),A(t)$ are nonnegative functions. The application of
these results to the geometric problems we shall consider below
leads us to assume the following structural conditions:
$$
\begin{array}{l}
\displaystyle A(t) \in L^\infty_\mathrm{loc}([0,+\infty)) \quad ,
\quad A(t)\ge
0 \quad , \quad A(t) \not\equiv 0 \\[0.2cm]
\displaystyle 0 \le v(t) \in L^\infty_\mathrm{loc}([0,+\infty))
\quad , \quad 1/v(t) \in L^\infty_\mathrm{loc}((0,+\infty))
\\[0.2cm]
\displaystyle v(t) \text{ is non decreasing near } 0 \text{ and }
\lim_{t\rightarrow 0^+} v(t) = 0.
\end{array}
$$
Of course, requests $A,v\ge 0$, $A\not\equiv 0$ are intended in
$L^\infty_\mathrm{loc}$ sense, while the last request means that
there exists a version of $v(t)$ which is non-decreasing in a
neighborhood of zero and whose limit as
$t\rightarrow 0^+$ is equal to zero.\\
Due to the weak regularity of $v$ and $A$, solutions $z(t)$ of \eqref{cauchy} are not expected to be classical,
and the Cauchy problem is expected to hold almost everywhere (a.e.) on $(0,+\infty)$. Equivalently (integrating
and using the condition in zero), we are interested in solutions $z(t)$ of the integral equation
$$
z(t) = z_0 - \int_0^t \frac{1}{v(s)} \left\{ \int_0^s A(x) v(x) z
(x) dx \right\} ds.
$$
For our purposes we shall look for $z(t)\in
\mathrm{Lip}_\mathrm{loc}([0,+\infty))$, that is, locally
Lipschitz solutions. Note that the locally Lipschitz condition
near zero ensures that $z'(t)=O(1)$ hold almost everywhere in a
neighborhood of zero. The existence of such solutions in our
assumptions will be given in Appendix, where we will also prove
that the zeros of $z(t)$, if any, are attained at isolated points.\\
The Cauchy problem \eqref{cauchy} is a somewhat "integrated" version of that presented in \cite{BR}, in the
sense that, as we shall see, in the geometric applications the role of $v(t)$ will be played by the volume
growth of geodesic spheres of some complete Riemannian manifold $M$, and $A(t)$ will represent the spherical
mean of some given function $a(x)$. However, the techniques introduced here are completely different from those
in \cite{BR}, and remind some in the work of Do Carmo and Zhou \cite{DoCZ}.\par
Nevertheless, as in \cite{BR}, we recognize an explicit
\emph{critical function} $\chi(t)$, depending only on $v(t)$,
which serves as a border line for the behavior of $z(t)$: roughly
speaking, if $A(t)$ is much greater than $\chi(t)$ in some region,
then $z(t)$ has a first zero, while if $A(t)$ is not greater than
$\chi(t)$ there are examples of positive solutions. We will see
that $\chi(t)$ generalizes the critical functions presented in
\cite{BR}.\par
Using $\chi(t)$ we will provide a condition in finite form for the
existence and localization of a first zero of $z(t)$ (Corollary
\ref{cor:zero}), and a sharp condition for the oscillatory
behavior (Corollary \ref{corforte}). In particular, this latter
Corollary improves on the application of the Hille-Nehari
oscillation theorem (see \cite{RRS}) to \eqref{cauchy}.
\par
The key technical result of the paper is Theorem \ref{teo:princ}
which, under very general assumptions, estimates the distance
between two consecutive zeros of an oscillatory solution of
(\ref{cauchy}): denoting with $T_1(\tau)< T_2(\tau)$ the first two
consecutive zeros of $z(t)$ after $t=\tau$, Theorem
\ref{teo:princ} states that
$$
T_2(\tau)-T_1(\tau) = O(\tau) \qquad \mathrm{as \ } \tau
\rightarrow +\infty.
$$
This result is achieved using a new but elementary technique which
highly improves on the application of Sturm's type arguments to
(\ref{cauchy}). Roughly speaking, the estimate will be obtained
performing a careful control on the level sets of the solution of
the Riccati equation associated to (\ref{cauchy}). Moreover, in
case
$$
v(t) \le f(t) = \Lambda \exp\{at^\alpha\log^\beta t\} \qquad
\Lambda,a,\alpha
>0, \ \beta \ge 0,
$$
we provide an upper estimate for
$$
\limsup_{\tau\rightarrow +\infty} \frac{T_2(\tau)}{\tau}
$$
with an explicit constant depending only on $\alpha$ and the
growth of $A(t)$ with respect to $\chi(t)$ (more precisely, with
respect to a critical curve $\chi_f(t)$ modelled on $f(t)$ instead
of $v(t)$).\par
There are several geometric applications of the above results; the
main idea is that (\ref{cauchy}) naturally appears in spectral
estimates. We will follow two slightly different ways. On the one
hand, we will provide an index estimate for Schr\"odinger type
operators $L=\Delta + a(x)$, while, on the other hand, we will
bound from above the growth of the spectral radius of the
Laplacian outside geodesic balls, even when the volume growth of
the manifold is faster than exponential. Applications naturally
arise in the setting of minimal hypersurfaces of Euclidean space,
their Gauss map, minimal surfaces and the Yamabe problem. We state
these geometric results in the next subsections.
\subsection{The geometric setting}
From now on, we let $(M,\langle,\rangle)$ denote a connected,
geodesically complete, non compact Riemannian manifold of
dimension $m\ge 2$. Fix an origin $o\in M$ and let $r(x)=
\mathrm{dist}(x,o)$ be the distance function from $o$. It is well
known that $r(x)$ is a Lipschitz function on $M$ which is smooth
outside $o$ and its cut-locus $\mathrm{cut}(o)$. For later use we
briefly recall some basic facts on the cut-locus in case $M$ is
geodesically complete; the interested reader
can consult, for instance, \cite{DoC} pp. 267-275.\\
Denote with $\exp$ the exponential map
$$
\exp \ : \ T_oM \rightarrow M,
$$
which, by the Hopf-Rinow theorem, is surjective and defined on the
whole $T_oM$. The origin $o$ is called a pole of $M$ if it has no
conjugate points; for example, this is the case if the sectional
curvature of $M$ is non positive. It turns out that, if $o$ is a
pole, $\exp$ is a covering map, hence a diffeomorphism if $M$ is
simply connected. For every $w \in T_oM$ such that $|w|=1$, we
indicate with $\gamma_w: [0,+\infty)\rightarrow M$ the geodesic
ray starting from $o$ with velocity $1$ in the direction of $w$,
and we consider
$$
t_w = \sup \left\{ s\in [0,+\infty) \ \text{ such that } \
r(\gamma_w(s)) = s \right\}.
$$
Clearly, $t_w>0$ because of the existence of geodesic
neighborhoods. If $t_w < +\infty$, we define the cut-point of $o$
along $\gamma_w$ as $\gamma_w(t_w)$. The cut-locus of $o$ is
defined as the union of the cut-points of $o$ along every geodesic
ray. In other words, $\mathrm{cut}(o)=\exp( \Sigma)$, where
$$
\Sigma = \left\{ tw \in T_oM \ : \ |w|=1 \ \text{ and } \ t = t_w
< +\infty \right\}.
$$
It is easy to see that, if $r(\gamma_w(s))=s$ for some $s>0$, then
the same equality holds for every $t\in [0,s)$. Therefore, if
$t_w<+\infty $ then $\gamma_w$ is length minimizing for every $t
\in (0,t_w]$ and it does not minimize length for any $t\in
(t_w,+\infty)$. By the Hopf-Rinow theorem we argue that the
exponential map restricted to the set $\mathcal{U} \cup \Sigma$,
where
$$
\mathcal{U} = \left\{ tw \in T_oM \ : \ |w|=1 \ \text{ and } \ t <
t_w \right\}
$$
is still surjective, hence $\exp(\mathcal{U})= M \backslash
\mathrm{cut}(o)= \mathrm{cut}(o)^c$. One can prove that
\begin{itemize}
\item[-] $\mathrm{cut}(o)$ is a zero measure, closed subset of
$M$, hence $\mathcal{U}=\exp^{-1}\big\{\mathrm{cut}(o)^c\big\}$ is
open in $T_oM$;
\item[-] $M$ is compact if and only if, for every $w \in T_oM$,
$t_w<+\infty$.
\item[-] $p\in \mathrm{cut}(o)$ if and only if either it is a
    conjugate point of $o$, or there exist at least 2 distinct
    geodesics joining $o$ to $p$ with the same length. The two possibilities do not
    reciprocally exclude.
\item[-] For every $q\in \exp(\mathcal{U})$, there exists a unique
minimizing geodesic from $o$ to $q$. In other words, $\exp :
\mathcal{U} \rightarrow \mathrm{cut}(o)^c$ is a bijection (indeed,
a diffeomorphism).
\end{itemize}
We indicate with $B_r$ the geodesic ball of radius $r$ centered at
$o$, with $\partial B_r$ its boundary and we call $\partial B_r
\cap \mathrm{cut}(o)^c$ the regular part of $\partial B_r$. The
regular part of $\partial B_r$ is an open set in the induced
topology on $\partial B_r$, and $\partial B_r \cap
\mathrm{cut}(o)^c$ is diffeomorphic, through the exponential map,
to the set $\mathcal{U} \cap \mathbb{S}^{m-1}(r)$, where
$\mathbb{S}^{m-1}(r)$ is the hypersphere
$$
\mathbb{S}^{m-1}(r) = \left\{ w \in T_oM \ : \ |w|=r \right\}.
$$
We denote with $\mathrm{Vol}(\partial B_r)$ the
$(m-1)$-dimensional volume of $\partial B_r$, that is, the
Hausdorff measure of $\partial B_r$. It turns out that it
coincides with the induced Riemannian measure when restricted to
the regular part of $\partial B_r$. The points of $\partial B_r
\cap \mathrm{cut}(o)$ may be image of many points of $\Sigma \cap
\mathbb{S}^{m-1}(r)$. For this reason, indicating with $\theta$ a
point of the unit sphere $\mathbb{S}^{m-1}=\mathbb{S}^{m-1}(1)
\subset T_oM$, we define the multiplicity function
$$
n_r(\theta) = \text{ cardinality of } \{ \varphi \in
\mathbb{S}^{m-1} \ : \ \exp(r\theta) = \exp(r\varphi) \} \le
+\infty.
$$
This coincides with the number of distinct minimizing geodesic
segments joining $o$ to $q=\exp(r\theta)$, which, analogously, we
denote with $n_r(q)$.
%
%
According to the work of Grimaldi and Pansu \cite{GP}, if we set
$$
\chi_r(\theta) = \left\{\begin{array}{ll} 1 & \text{if } \ r <
t_\theta \\[0.1cm]
1/n_r(\theta) & \text{if } \ r =
t_\theta \\[0.1cm]
0 & \text{if } \ r > t_\theta
\end{array} \right. \qquad \text{and} \qquad \chi_{r^\pm}(\theta) =
\lim_{t\rightarrow r^\pm} \chi_t(\theta)
$$
then the Hausdorff measure of $\partial B_r$ is given by
$$
\mathrm{Vol}(\partial B_r) = \int_{\mathbb{S}^{m-1}}
\Theta(r,\theta) \chi_r(\theta) d\theta,
$$
where $\Theta(r,\theta)$ is the density of the Riemannian measure.
Moreover, by the dominated convergence theorem,
$$
\lim_{t\rightarrow r^\pm}\mathrm{Vol}(\partial B_t) =
\int_{\mathbb{S}^{m-1}} \Theta(r,\theta) \chi_{r^\pm}(\theta)
d\theta.
$$
Therefore, in general circumstances $\mathrm{Vol}(\partial B_r)$
may present discontinuities of the ``first kind", that is, at a
point $r>0$ we always have the existence of finite limits both
from the right and from the left, possibly with two different
values. Indeed, setting $v(r)=\mathrm{Vol}(\partial B_r)$, it is
shown in \cite{GP} that for every complete Riemannian manifold
$$
v(r^+) - v(r^-) = -2 \mathrm{Vol}(\partial B_r \cap
\mathrm{cut}(o)).
$$
The key ingredient of their proof is a technical Lemma which shows
that, up to a set of $(m-1)$-dimensional measure zero, $\partial
B_r \cap \mathrm{cut}(o)$ is made up of points having exactly $2$
distinct geodesics which minimize distance from $o$. Observe that
$v(t)$ jumps downward and that, a priori, the discontinuities of
$v(t)$ may be non isolated. Note also that, from the definition of
$\chi_r(\theta)$, we get
$$
\chi_{r^+} (\theta) = \left\{\begin{array}{ll} 1 & \text{if } \ r
<
t_\theta \\[0.1cm]
0 & \text{if } \ r =
t_\theta \\[0.1cm]
0 & \text{if } \ r > t_\theta
\end{array} \right. \qquad
\chi_{r^-}(\theta) = \left\{\begin{array}{ll} 1 & \text{if } \ r <
t_\theta \\[0.1cm]
1 & \text{if } \ r =
t_\theta \\[0.1cm]
0 & \text{if } \ r > t_\theta
\end{array} \right. ,
$$
hence from $\chi_{r^+} \le \chi_r \le \chi_{r^-}$ we deduce that
$v(t) \in [v^+(t),v^-(t)]$. Therefore, a necessary and sufficient
condition on $\mathrm{Vol}(\partial B_r)$ to be continuous on
$[0,+\infty)$ is given by the ``transversality condition"
$$
\mathrm{Vol}(\partial B_r \cap \mathrm{cut}(o)) = 0 \qquad \forall
\ r \ge 0.
$$
However, this reasonable request sometimes is not easy to verify.
This is the case, for example, when one constructs manifolds as
immersed submanifolds of some ambient space. This suggests to work
with discontinuous volume functions $\mathrm{Vol}(\partial B_r)$
which will take the role of $v$ in \eqref{cauchy}. The next result
will reveal important in what follows.
\begin{Proposizione}\label{proprietavol}
Let $v(r) = \mathrm{Vol}(\partial B_r)$ be the volume of geodesic
spheres of a connected, complete, non compact Riemannian manifold.
Then $v(r)$ is continuous and increasing in a neighborhood of
$r=0$. Furthermore,
\begin{equation}\label{proprvol}
v(r) = \frac{v(r^+) + v(r^-)}{2} \quad , \quad v(r)
> 0 \ \text{ for } \ r>0 \quad , \quad \frac{1}{v(r)} \in
L^\infty_{\mathrm{loc}}((0,+\infty))
\end{equation}
\end{Proposizione}
\noindent \textbf{Proof. } The first part is immediate using polar
coordinates around zero. As for the first property in
\eqref{proprvol}, we denote with $\mathcal{V} = \{ w \in
\mathbb{S}^{m-1} \ : \ rw \in \mathcal{U} \}$ and with
$\mathcal{W} = \{ w \in \mathbb{S}^{m-1} \ : \ rw \in \Sigma\}$.
Since $r\mathcal{V} = \mathbb{S}^{m-1}(r) \cap \mathcal{U}$ is
open, then $\mathcal{V}$ is an open set of $\mathbb{S}^{m-1}$. In
polar coordinates
$$
\begin{array}{ll}
v(r) & = \displaystyle \int_{\mathbb{S}^{m-1}} \Theta(r,\theta)
\chi_r(\theta) d\theta \equiv \int_{\mathcal{V}} \Theta(r,\theta)
d\theta + \int_{\mathcal{W}}
\Theta(r,\theta) \frac{1}{n_r(\theta)} d\theta \\[0.5cm]
 & = \mathrm{Vol}(\partial B_r \cap \mathrm{cut}(o)^c) +
 \mathrm{Vol}(\partial B_r \cap \mathrm{cut}(o)) \\[0.5cm]
v(r^+) & = \displaystyle \int_{\mathbb{S}^{m-1}} \Theta(r,\theta)
\chi_{r^+}(\theta) d\theta \equiv \int_{\mathcal{V}}
\Theta(r,\theta)
d\theta = \mathrm{Vol}(\partial B_r \cap \mathrm{cut}(o)^c) \\[0.5cm]
v(r^-) & =  \displaystyle \int_{\mathbb{S}^{m-1}} \Theta(r,\theta)
\chi_{r^-}(\theta) d\theta \equiv \int_{\mathcal{V}}
\Theta(r,\theta)
d\theta + \int_{\mathcal{W}} \Theta(r,\theta) d\theta \\[0.5cm]
& = \displaystyle \mathrm{Vol}(\partial B_r \cap
\mathrm{cut}(o)^c) + \int_{\partial B_r \cap \mathrm{cut}(o)}
n_r(x)d\sigma(x)
\end{array}
$$
By the Grimaldi-Pansu lemma \cite{GP}, up to a set of
$(m-1)$-dimensional measure zero, the multiplicity $n_r(x)$ is
equal to $2$. Therefore, by the above expressions is immediate to
deduce that $v(r^+)+v(r^-) = 2 v(r)$. We observe now that if we
prove that $1/v \in L^\infty_{\mathrm{loc}}((0,+\infty))$, then
$v(r)>0$ on $(0,+\infty)$. Indeed, assume $v(r_0)=0$ for some
$r_0\in (0,+\infty)$. Then necessarily $v(r_0^+)=0$, $v(r_0^-) =
2v(r_0)-v(r_0^+)=0$ and $1/v$ is unbounded in a neighborhood of
$r_0$. It remains to prove that $1/v \in
L^\infty_{\mathrm{loc}}((0,+\infty))$, that is, $v(r)$ is bounded
away from zero on every compact set $K$ disjoint from $r=0$.
Assume by contradiction that there exists $\{r_k\} \subset K$ such
that $v(r_k) \rightarrow 0$. By compactness, there exists
$\widetilde{r}\in K$ such that $r_k \rightarrow \widetilde{r}$. Up
to passing to a subsequence we have two cases: $r_k \uparrow
\widetilde{r}$ or $r_k \downarrow \widetilde{r}$. In the first
case $v(\widetilde{r}^-)=0$, in the second $v(\widetilde{r}^+)=0$.
However, since $v$ jumps downward, in both cases
$v(\widetilde{r}^+)=0$. We are going to show that
\begin{equation}\label{absu}
\partial B_{\widetilde{r}} \subseteq \mathrm{cut}(o).
\end{equation}
Indeed, let \eqref{absu} be false, and let $q\in \partial
B_{\widetilde{r}} \cap \mathrm{cut}(o)^c$. Since $\exp$ is a
diffeomorphism in a neighborhood of $q$, we can choose a unique
$\theta_0 \in \mathcal{V}$ such that $q =
\exp(\widetilde{r}\theta_0)$. Moreover, since $\mathcal{U}$ is
open, from $\widetilde{r}\theta_0 \in \mathcal{U}$ we can chose a
neighborhood $\mathcal{J}$ with compact closure in $\mathcal{U}$
of the form
$$
\mathcal{J} = \{ r\theta \ : \ r \in
(\widetilde{r}-2\varepsilon,\widetilde{r}+2\varepsilon) \ , \
\theta \in \mathcal{V}_{\theta_0} \},
$$
where $\varepsilon>0$ is sufficiently small and
$\mathcal{V}_{\theta_0}$ is a neighborhood of $\theta_0$ on the
unit sphere $\mathbb{S}^{m-1}$, independent from $\varepsilon$.
Since the Riemannian density $\Theta$ is smooth and positive,
there exists $C>0$ independent of $\varepsilon$ such that
$\Theta(r,\theta) \ge C$ on $\mathcal{J}$. It follows that
$$
v(\widetilde{r}+\varepsilon) = \int_{\mathbb{S}^{m-1}}
\Theta(\widetilde{r}+\varepsilon,\theta)
\chi_{\widetilde{r}+\varepsilon}(\theta) d\theta \ge
\int_{\mathcal{V}_{\theta_0}} C d\theta = C
\mathrm{Vol}_{\mathrm{Eucl}}(\mathcal{V}_{\theta_0}) \qquad
\forall \ \varepsilon.
$$
This contradicts $v(\widetilde{r}^+)=0$ and proves \eqref{absu}.
By \eqref{absu} we deduce that, for every geodesic ray $\gamma_w$
starting from $o$, there exists $t_w\le r$ such that
$\gamma_w(t_w) \in \mathrm{cut}(o)$. Therefore, $M$ is compact
with diameter $\le 2r$, against our assumptions. $\blacktriangle$\\
\par
Let $s(x)$ be the scalar curvature of $(M,\langle,\rangle)$. The
previous proposition enables us to define the spherical mean
$$
S(r) = \frac{1}{\mathrm{Vol}(\partial B_r)} \int_{\partial B_r} s
$$
on the whole $(0,+\infty)$. $S(r)$ is continuous in a neighborhood
of zero with $\lim_{r\rightarrow 0^+} S(r) = s(o)$, and possesses
at least the same regularity as $\mathrm{Vol}(\partial B_r)$. In
case $(\mathrm{Vol}(\partial B_r))^{-1} \in L^1(+\infty)$ we
define the \textbf{critical function}
\begin{equation}\label{defchi}
\chi(r) = \Big( 2 \mathrm{Vol}(\partial B_r) \int_r^{+\infty}
\frac{ds}{\mathrm{Vol}(\partial B_s)} \Big)^{-2} \ \ \in
L^\infty_{\mathrm{loc}}((0,+\infty)),
\end{equation}
that we shall consider below.

Since in the sequel we will be concerned with spectral arguments,
we briefly recall some definitions. Let $\Delta$ denote the
Laplace-Beltrami operator on $M$, and consider a differential
operator $L= \Delta + a(x)$, where $a(x) \in C^0(M)$, and a
bounded domain $\Omega \subset M$. The $k$-th eigenvalue
$\lambda_k^L (\Omega)$, of $L$ on $\Omega$ (counted with its
multiplicity) is defined by Rayleigh characterization:
\begin{equation}\label{eq:uno}
\lambda_k^L (\Omega)= \inf_{ \scriptsize{\begin{array}{c} V_{k}
\le C^\infty_0(\Omega) \\[0.1cm] \dim(V_{k})=k \end{array} } } \left(
\sup_{ 0 \neq \phi \in V_{k}}  \frac{\int_\Omega {| \nabla \phi
|^2} - \int_\Omega a\phi^2}{\int_\Omega \phi^2}\right),
\end{equation}
%
%
where we can substitute $C^\infty_0(\Omega)$ with
$\mathrm{Lip}_0(\Omega)$. If $\Omega$ has sufficiently regular
boundary, $\lambda_1^L (\Omega)$ is achieved by the non zero
solutions of the Dirichlet problem
\begin{equation}\label{eq:due}
\left\{
\begin{array}{l}
Lu + \lambda_1^L (\Omega) u= 0 \qquad {\rm on} \ \Omega \\
u\equiv 0 \qquad {\rm on} \ \partial \Omega
\end{array}
\right.
\end{equation}
Note that $L$ is non positive on $C^\infty_0(\Omega)$ if and only
if $\lambda_1^L(\Omega) \ge 0$. The main example of a non positive
operator on every $\Omega$ is the Laplacian itself.\par
We define the index $\mathrm{ind}_L(\Omega)$ as the number of
negative eigenvalues of $-L$. By Rellich Theorem, this number is
finite. Indeed, using Rayleigh characterization
$$
\lambda_k^L(\Omega) \ge \lambda_k^\Delta(\Omega) - \|a\|_{L^\infty(\Omega)},
$$
therefore $\widetilde{L}= L - \|a\|_{L^\infty(\Omega)}$ is
strictly non positive on $C^\infty_0(\Omega)$, hence it is
invertible. The Friedrich extension of $(-\widetilde{L})^{-1}:
L^2(\Omega) \rightarrow L^2(\Omega)$ is a compact operator, so
that its spectrum consists in a discrete sequence $\{\lambda_j\}$
of eigenvalues, each of them with finite multiplicity. It follows
that the spectrum of $-L$ is $\{\lambda_j-\|a\|_{L^\infty(\Omega)}
\}$, and $\mathrm{ind}_L(\Omega)$ is clearly finite.
The bottom of the spectrum of $L$ on $M$, also called the first
eigenvalue or the spectral radius, $\lambda_1^L (M)$, is defined
by
\begin{equation}\label{eq:tre}
\lambda_1^L (M)=\inf \left\{ \lambda_1^L (\Omega) \ : \ \Omega
\subset M \ {\rm is} \ {\rm a } \ {\rm bounded} \ {\rm domain}
\right\}
\end{equation}
Let $Z \subset M$ be a subset. We define the first eigenvalue of
$L$ on the "punctured" manifold $M\backslash Z$ by
\begin{equation} \label{eq:quattro}
\lambda_1^L (M\backslash Z)=\inf \left\{ \lambda_1^L (\Omega) \ :
\ \Omega \subset M\backslash Z \ {\rm is} \ {\rm a } \ {\rm
bounded} \ {\rm domain} \right\}
\end{equation}
Similarly, the index of $L$ on $M$ is defined by
$$
\mathrm{ind}_L(M) = \sup \left\{ \mathrm{ind}_L (\Omega) \ : \
\Omega \subset M \ {\rm is} \ {\rm a } \ {\rm bounded} \ {\rm
domain} \right\}
$$
and it may be infinite. Note that $\mathrm{ind}_L(M) = 0$ if and
only if $\lambda_1^L(M) \ge 0$.
\subsection{Spectral estimates: the two main results}
The first theorem deals with the index of $L$.
\begin{Teorema}\label{spettro}
Let $a(x) \in C^0(M)$. Suppose that the spherical mean $A(r)$ of
$a(x)$ is non negative and not identically null. Consider the
following assumptions:
\begin{itemize}
\item[$(i)$] either
$$
(\mathrm{Vol}(\partial B_r))^{-1}\not \in    L^1(+\infty)
$$
or $(\mathrm{Vol}(\partial B_r))^{-1} \in
    L^1(+\infty)$ and there exist $0<R_0<R_1$ such that $A(r)
    \not \equiv 0$ on $[0,R_0]$ and
\begin{equation}\label{finite}
\displaystyle \int^{R_1}_{R_0}
\big(\sqrt{A(s)}-\sqrt{\chi(s)}\big) ds  > -\frac{1}{2} \Big(\log
\int_{B_{R_0}} a +\log \int_{R_0}^{+\infty}
{\frac{ds}{\mathrm{Vol}(\partial B_s)}}\Big)
\end{equation}
\item[$(ii)$]  either
\begin{equation}
(\mathrm{Vol}(\partial B_r))^{-1} \not\in L^1(+\infty) \quad ,
\quad a(x) \not \in L^1(M)
\end{equation}
or
\begin{equation}\label{infinite}
(\mathrm{Vol}(\partial B_r))^{-1} \in L^1(+\infty) \quad , \quad
\limsup_{r\rightarrow +\infty} \int_R^r
\big(\sqrt{A(s)}-\sqrt{\chi(s)}\big) ds = +\infty.
\end{equation}
for some $R$ sufficiently large.
%
%
%
%
\item[$(iii)$] $(\mathrm{Vol}(\partial B_r))^{-1} \in
    L^1(+\infty)$,
$$
\mathrm{Vol}(\partial B_r) \le \Lambda \exp\{ar^\alpha \log^\beta
r\} \qquad \mathit{ for \ some \ } \Lambda,a,\alpha >0, \ \beta
\ge 0
$$
and, for some $R>0$, $c>1$,
\begin{equation}\label{growthA}
\sqrt{A(r)} \ge c \Big(\frac{a \alpha}{2}\Big) r^{\alpha
-1}\log^\beta r \qquad \forall \ r \ge R
\end{equation}
\end{itemize}
Let $L= \Delta + a(x)$. Then
\begin{itemize}
\item[-] under assumption $(i)$, $\lambda_1^L(M)<0$;
\item[-] under assumption $(ii)$, $L$ is unstable at infinity,
    that is, $\lambda_1^L(M\backslash B_R) < 0$ for every
    $R>0$. In particular, $L$ has infinite index;
\item[-] under assumption $(iii)$, $L$ is unstable at infinity
    and
\begin{equation}\label{indexgrowth}
\liminf_{r\rightarrow +\infty} \frac{\mathrm{ind}_L(B_r)}{\log r}
\ge \frac{\alpha}{2\log \Big(\frac{c+1}{c-1}\Big)}.
\end{equation}
\end{itemize}
\end{Teorema}
We observe that (\ref{infinite}) and (\ref{growthA}) are
conditions ``at infinity" and they are typical of oscillation
results. On the other hand, condition (\ref{finite}) deserves some
special attention since it is in finite form, in the sense that it
only involves the behavior of $a(x)$ on a compact set, namely
$B_{R_1}$: the left hand side states how much must $a(x)$ exceed
the critical curve on the compact annular region
$\overline{B}_{R_1}\backslash B_{R_0}$ in order to have a negative
spectral radius, and it only depends on the behavior of $a(x)$
near zero (on $B_{R_0}$) and on the geometry at infinity of $M$.
Note also that $R_1$ does not appear in the right
hand side of (\ref{finite}).
\begin{Remark}\label{indicestabile}
\emph{By a famous result of Fisher-Colbrie \cite{FC}, condition
$\mathrm{Ind}_L(M)<\infty$ implies the stability at infinity (that
is, $\lambda_1^L(M\backslash B_R)\ge 0$ for some $R \ge 0$). As
far as we know, it is yet an open problem to prove the converse,
or to provide an explicit counterexample. However, we remark that
a sufficient condition to have finite index is that the strict
inequality $\lambda_1^L(M\backslash B_R)>0$ hold for some $R$. For
a detailed account of spectral theory for Schr\"odinger operators
on Riemannian manifolds we refer the reader to \cite{PRS1}.}
\end{Remark}
\par
The second result can be probably regarded as the core of the
paper: it provides a sharp upper bound for the growth of
$\lambda_1^\Delta (M\backslash B_R)$ as a (monotone) function of
$R$. In the literature, bounds for the spectral radius on $M$ are
obtained under at most exponential volume growth of geodesic
spheres. On the contrary, Theorem \ref{principale} works also with
faster volume growths. To better appreciate the result that we
shall introduce below, we begin with some preliminary
considerations.\par
It is well known that, if $Z$ is any compact subset of $\R^m$,
then $\lambda_1^\Delta (\R^m\backslash Z)=0$. Extending a result
of Cheng and Yau \cite{CY}, Brooks \cite{B} has shown that if the
manifold $(M,\langle,\rangle)$ has at most sub-exponential volume
growth then $\lambda_1^\Delta (M)=0$. However, if we puncture the
manifold by a compact set $Z \not= \emptyset$, contrary to the
case of $\R^m$, it may happen that $\lambda_1^\Delta (M\backslash
Z) \not=0$. Indeed, Do Carmo and Zhou, \cite{DoCZ}, give an
example where $\mathrm{Vol} (M)<+\infty$ and
$$
\lambda_1^\Delta (M\backslash \overline{B}_1) \ge \frac{1}{4}
$$
Moreover, up to the missing requirement of continuity of
$\mathrm{Vol}(\partial B_r)$, they prove that \emph{in case $M$
has infinite volume},
\begin{itemize}
\item[-] if $M$ has sub-exponential volume growth of geodesic
spheres, then
\begin{equation}\label{stimedocz1}
\lambda_1^\Delta(M\backslash B_R) =0 \qquad \forall \ R \ge 0;
\end{equation}
\item[-] if $\mathrm{Vol}(\partial B_r) \le C{\rm e}^{ar}$ for
some $C,a>0$, then
\begin{equation}\label{stimedocz2}
\lambda_1^\Delta(M\backslash B_R)\le \frac{a^2}{4} \qquad \forall
\ R \ge 0.
\end{equation}
\end{itemize}
It is interesting to see what happens when the volume growth is
faster than exponential. Towards this aim, we extend Do Carmo and
Zhou's example to grasp the situation a step further. Thus we
consider the model, in the sense of Greene and Wu,
$(M,ds^2)=(\R^m, ds^2)$, with metric given in polar coordinates by
\begin{equation}\label{eq:cinque}
ds^2=dr^2+h(r)^2d\theta^2
\end{equation}
where $h \in C^\infty ([0,+\infty))$ is positive on $(0,+\infty)$
and satisfies
\begin{equation} \label{eq:sei}
h(r)= \left\{
\begin{array}{l}
\displaystyle r \quad {\rm on} \ [0,1]\\[0.2cm]
\displaystyle \exp\left\{\frac{ar^\alpha}{m-1}\right\} \quad {\rm
on} \ [2,+\infty)
\end{array}
\right.
\end{equation}
for some $a > 0, \ \alpha \ge 1$. Note that (\ref{eq:cinque})
extends smoothly at the origin because of the definition of $h$
near 0, and that, for $r\ge 2$, $\mathrm{Vol}(\partial B_r) =
\exp\{ar^\alpha\}$. We let $b\in (0,a)$ and set
\begin{equation} \label{eq:sette}
u_b(x)={\rm e}^{-b r(x)^\alpha} \qquad {\rm on} \ M\backslash B_2.
\end{equation}
A simple checking shows that
$$
\Delta u_b + \lambda_b (r) u_b = 0 \qquad {\rm on} \ M\backslash
B_2,
$$
where $\lambda_b(r)$ is defined as
\begin{equation}\label{eq:nove}
\lambda_b (r)= \alpha^2 b (a-b) r^{2(\alpha-1)}
+\alpha(\alpha-1)br^{\alpha-2}.
\end{equation}
Observe that, in case $\alpha=1$, $\lambda_b(r)\equiv b(a-b)$, while, if $\alpha>1$, $\lambda_b(r)$ is strictly
increasing on $(R_0,+\infty)$, with $R_0$ sufficiently large that
$$
2\alpha (a-b) R_0^\alpha + (\alpha-2) > 0.
$$
Up to further enlarging $R_0$, we can also assume that
\begin{equation} \label{eq:otto}
\frac{\alpha -1}{2 \alpha} \frac{1}{r^\alpha} < \frac{a}{2} \qquad
{\rm for} \ r \ge R_0.
\end{equation}
%
%
Applying a result of Cheng and Yau, \cite{CY} we have that, for
every $b \in (0,a)$, $R\ge R_0$,
$$
\lambda_1^\Delta (M\backslash B_R)\ge \inf_{M \backslash B_R}
{-\frac{\Delta u_b}{u_b}}=\inf_{[R,+\infty)} {\lambda_b (r)}
=\lambda_b (R)
$$
The choice
$$
\widetilde{b}=\frac{a}{2}+\frac{\alpha-1}{2
\alpha}\frac{1}{R^\alpha}
$$
maximize $\lambda_b (R)$ and $\widetilde{b}\in (0,a)$ because of
\eqref{eq:otto}. Then, for $R \ge R_0$
\begin{equation} \label{eq:dieci}
\lambda_1^\Delta (M \backslash B_R) \ge \alpha^2 \left(
\frac{a^2}{4}-\frac{(\alpha-1)^2}{4\alpha^2} \frac{1}{R^{2\alpha}}
\right) R^{2(\alpha-1)}
\end{equation}
Note that for $\alpha=1$ the above reduces to
$$
\lambda_1^\Delta (M \backslash B_R) \ge \frac{a^2}{4}
$$
In particular, this shows that the upper bound in Theorem 3.1 in
\cite{DoCZ} is sharp. This example, for $\mathrm{Vol}(\partial
B_r) \le C{\rm e}^{ar^\alpha}$, $C,a>0$, $\alpha\ge 1$, suggests
to look for an upper bound of $\lambda_1^\Delta (M \backslash
B_R)$ of the form
$$
C_1 R^{2 (\alpha-1)}
$$
with $C_1=C_1 (a,\alpha)>0$. The guess is indeed correct, as
Theorem \ref{principale} shows.
\begin{Teorema} \label{principale}
If $M$ is a connected, complete, non compact Riemannian manifold
such that
$$
(\mathrm{Vol}(\partial B_r))^{-1} \in L^1(+\infty) \quad , \quad \mathrm{Vol} (\partial B_r) \le \Lambda \exp[ar^\alpha \log^\beta
r] \qquad {\rm for} \ r \ {\rm large,}
$$
for some $\Lambda, a,\alpha >0, \ \beta \ge 0$, the following
estimates hold:
\begin{itemize}
\item[-] If $0<\alpha<1$ then
$$
\lambda_1^\Delta(M\backslash B_R)= 0 \qquad \forall \ R \ge 0.
$$
\item[-] If $\alpha \ge 1$ then
\begin{equation}\label{tesiuno}
\limsup_{R\to +\infty} \Bigg(\frac{\lambda_1^\Delta(M\backslash
B_R)}{R^{2(\alpha-1)}\log^{2\beta}R} \Bigg)\le
\frac{a^2\alpha^2}{4} \inf_{c\in (1,+\infty)} \Big\{
c^2\Big(\frac{c+1}{c-1}\Big)^{\frac{4(\alpha-1)}{\alpha}} \Big\}
\end{equation}
\end{itemize}

\end{Teorema}

\begin{Remark}\label{inL1}
\emph{Note that $(\mathrm{Vol}(\partial B_r))^{-1} \in
L^1(+\infty)$ implies $\mathrm{Vol}(M) = \infty$. This follows
from Schwarz inequality
$$
\int_r^R \frac{ds}{\mathrm{Vol}(\partial B_s)} \int_r^R
\mathrm{Vol}(\partial B_s)ds \ge (R-r)^2
$$
letting $R\rightarrow +\infty$.}
\end{Remark}
We stress that the hypothesis $\mathrm{Vol}(M)=\infty$ is
essential. In fact, Do Carmo and Zhou example quoted above shows
that the theorem fails if $\mathrm{Vol}(M) <\infty$. On the
contrary, the stronger assumption $(\mathrm{Vol}(
\partial B_r ))^{-1} \in L^1 (+\infty)$ is for convenience: if it
fails, we will show in Lemma \ref{cor:cinque} that $\lambda_1^\Delta (M \backslash B_R) =0$ for every $ R \ge
0$. We underline that in Theorem \ref{principale} we have been considering volume growth assumptions, which are
weaker and more general than the usual curvature conditions used in estimating $\lambda_1^\Delta(M)$ (see for
instance \cite{pin}). It is also worth mentioning that the problem of estimating $\lambda_1^\Delta(M\backslash
B_R)$ from above arises naturally in the study of unstable hypersurfaces with constant mean curvature: see for
example \cite{DoCZ} for details and further references.
\subsection{Geometric consequences}
The first geometric consequence is the following density theorem
for complete minimally immersed hypersurfaces of Euclidean space.
\begin{Teorema}\label{main}
Let $\varphi: M \rightarrow \mathbb{R}^{m+1}$ be a minimal
hypersurface. We identify $T_xM$ with $\varphi_\ast T_xM$ viewed
as an affine hyperplane in $\mathbb{R}^{m+1}$ passing through
$\varphi(x)$. Assume that
\begin{equation}\label{volumenoninl1}
(\mathrm{Vol}(\partial B_r))^{-1} \not\in L^1(+\infty) \quad , \quad s(x) \not\in L^1(M)
\end{equation}
or that
\begin{equation}\label{volume}
(\mathrm{Vol}(\partial B_r))^{-1} \in L^1(+\infty)  \quad , \quad \mathrm{Vol}(\partial B_r) \le \Lambda\exp\{r^\alpha\}\\
\end{equation}
\begin{equation}\label{mediasferica}
S(r) \le -\frac{C}{r^\mu}
\end{equation}
for $r\gg 1$ and some constants $C,\Lambda,\alpha>0$, $\mu\in
\mathbb{R}$, with
\begin{equation}\label{parametri}
2\alpha < 2-\mu.
\end{equation}
Then, for every compact set $\Omega\subseteq M$
\begin{equation}\label{tuttoRm}
\bigcup_{x\in M\backslash \Omega} T_xM \equiv \mathbb{R}^{m+1}
\end{equation}
\end{Teorema}
We note that Halpern \cite{H} has proved that, when the
hypersurface is compact and orientable, $\bigcup_{x \in M} T_xM
\not\equiv \mathbb{R}^{m+1}$ if and only if $M$ is embedded as the
boundary of an open star-shaped domain of $\mathbb{R}^{m+1}$. In
case $M$ is non compact there are many examples with $\bigcup_{x
\in M} T_xM \not\equiv \mathbb{R}^{m+1}$, for instance cylinders
over suitable curves. However, in case $m=2$ complete minimal
surfaces in $\mathbb{R}^3$ for which $\bigcup_{x\in M} T_xM
\not\equiv \mathbb{R}^3$ are planes: this has been proved
by Hasanis and Koutroufiotis in \cite{HK}.\\
\par
In an analogous way, we prove the following result. 
%
\begin{Teorema} \label{Gaussmap}
Let $\varphi : M \rightarrow \mathbb{R}^{m+1}$ be a connected,
complete non compact minimal hypersurface. Assume that either
\begin{equation}\label{volumenoninl12}
(\mathrm{Vol}(\partial B_r))^{-1} \not\in L^1(+\infty) \quad ,
\quad s(x) \not\in L^1(M)
\end{equation}
or that, for some $C,\Lambda,\alpha>0$, $\mu \in \mathbb{R}$
\begin{equation}\label{volescalare}
\begin{array}{l}
\displaystyle (\mathrm{Vol}(\partial B_r))^{-1} \in L^1(+\infty)
\quad , \quad \mathrm{Vol}(\partial B_r) \le
\Lambda\exp\{r^\alpha\}\\[0.3cm]
\displaystyle S(r) \le -\frac{C}{r^\mu} \qquad \text{and} \qquad
2\alpha<2-\mu.
\end{array}
\end{equation}
Fix an equator $E$ in $\mathbb{S}^m$. Then the spherical Gauss map
$\nu$ meets $E$ infinitely many times along a divergent sequence
in M.
\end{Teorema}
Note that we have not assumed the orientability of $M$; hence, the
spherical Gauss map is only locally defined. However, due to the
central symmetry of the equators, the conclusion of the theorem
does not depend on the chosen local
orientation: if $\nu(x) \in E$, then also $-\nu(x) \in E$.\\
\par
As a third consequence of Theorem \ref{spettro}, we have the
following result of Fisher-Colbrie \cite{FC} and Gulliver
\cite{G}.
\begin{Teorema}\label{superfici}
%
%
Let $N$ be a flat $3$-manifold, and let $\varphi : M \rightarrow
N$ be a simply connected, minimally immersed surface. We denote with $K$ the
(necessarily non positive) sectional curvature of $M$.\\
Consider the stability operator $L=\Delta + |II|^2$. If $M$ is
stable at infinity (in particular, if $\mathrm{Ind}_L(M)<
\infty$), then $M$ is parabolic and
\begin{equation}\label{supminimali}
\int_M |K| < +\infty.
\end{equation}
\end{Teorema}
%
%
\begin{Remark}
\emph{Indeed, it is shown in \cite{FC} that
$$
\int_M |K| < +\infty \Longleftrightarrow \mathrm{Ind}_L(M) <
\infty.
$$
Combining with the above result yields that if a simply connected, minimal surface in $\mathbb{R}^3$ is stable
at infinity, then its stability index is finite. Therefore, referring to Remark \ref{indicestabile}, a possible
counterexample to the equivalence between $\mathrm{Ind}_L(M)<\infty$ and $\lambda_1^L(M\backslash B_R) \ge 0$
must not be searched in the setting of
 minimal surfaces in $\mathbb{R}^3$, when $L$ is the stability
operator.}
\end{Remark}
With the same technique, we recover a well known result of Do
Carmo and Peng \cite{DoCP}, Fisher-Colbrie and Schoen \cite{FCS}
and Pogorelov \cite{Po}.
\begin{Corollario}\label{corpiatto} Let $\varphi: M\rightarrow
\mathbb{R}^3$ be a minimally immersed surface. If $M$ is stable,
then $M$ is totally geodesic (hence, an affine plane).
\end{Corollario}
%
%
%
%
%
%
%
\par
The last geometrical application employs directly Theorem
\ref{spettro}, together with Theorems $2.4$ and $2.1$ of
\cite{PRS}, to yield the following existence result for the Yamabe
problem which requires no assumptions on the Ricci curvature.
\begin{Teorema}\label{Yamabe}
Suppose that the dimension of $M$ is $m\ge 3$ and that the
spherical mean $S(r)$ satisfies
$$
S(r) \le 0 \quad \mathit{on \ } [0,+\infty) \quad ,
\quad S \not \equiv 0
$$
Let $k(x)\in C^\infty(M)$ be non positive on $M$ and strictly
negative outside a compact set. Set $\mathcal{K}_0 = k^{-1}\{0\}$
and, for
$$
L= \Delta - \frac{1}{c_m}s(x)  \qquad \mathit{where} \quad c_m =
\frac{4(m-1)}{m-2},
$$
define $\lambda_1^L(\mathcal{K}_0) = \sup_D \lambda_1^L(D)$, where
$D$ varies among all open sets with smooth boundary containing
$\mathcal{K}_0$. Suppose
$$
\lambda_1^L(\mathcal{K}_0) >0.
$$
Assume that either $(\mathrm{Vol}(\partial B_r))^{-1}\not \in
L^1(+\infty)$ or otherwise that there exists $0< R_0 < R_1$ such
that $S\not\equiv 0$ on $[0,R_0]$ and
\begin{equation}\label{yam}
\displaystyle \int^{R_1}_{R_0}
\big(\sqrt{\frac{|S(t)|}{c_m}}-\sqrt{\chi(t)}\big) dt
> -\frac{1}{2} \Big(\log \int_{B_{R_0}} \frac{|s(x)|}{c_m} +\log
\int_{R_0}^{+\infty} {\frac{dt}{\mathrm{Vol}(\partial B_t)}}\Big)
\end{equation}
Then, the metric $\langle,\rangle$ can be conformally deformed to
a
new metric of scalar curvature $k(x)$.\\
\end{Teorema}
As the discussion after Theorem \ref{spettro} suggests, this
latter result implies that a strongly negative scalar curvature on
a compact region $\Omega$ gives the existence of the conformal
deformation independently of the behavior of $s(x)$ outside
$\Omega$.
%
%
\section{Existence of a first zero and oscillations}\label{existencefirstzero}
Fix $R\in (0,+\infty]$ (note that the value $+\infty$ is allowed),
and consider the following set of assumptions:
\begin{align}
&0 \le A(t) \in L^\infty_{\mathrm{loc}}([0,R))
\quad , \quad A\not\equiv 0 \quad \text{in } L^\infty_\mathrm{loc} \text{ sense} & & \tag{A1} \label{A}     \\
&0 \le v(t) \in L^\infty_{\mathrm{loc}}([0,R)) \quad, \quad
\frac{1}{v(t)} \in L^\infty_{\mathrm{loc}}((0,R)) \quad , \quad
\lim_{t\rightarrow 0^+} v(t)=0 & & \tag{V1} \label{V}
\end{align}
In case $1/v \in L^1(R^-)$, we define the critical function
\begin{equation}\label{chi}
\chi_R(t) = \Big(\displaystyle 2v(t)\int_t^{R}
\frac{ds}{v(s)}\Big)^{-2} = \Big[\Big( -\frac{1}{2} \log\int_t^{R}
\frac{ds}{v(s)}\Big)'\Big]^2 \quad \in
L^\infty_\mathrm{loc}((0,R))
\end{equation}
For the ease of notation we write $\chi(t)$ in case $R=+\infty$.
We are now ready to prove:
\begin{Teorema} \label{teo:uno}
Let $A,v$ satisfy \eqref{A}, \eqref{V} and let $z \in
\mathrm{Lip}_\mathrm{loc}([0,R))$ be a positive solution of
\begin{equation}\label{cauchy2}
\left\{ \begin{array}{l} (v(t)z'(t))' + A(t)v(t) z(t) = 0 \qquad
\text{almost everywhere on } (0,R)\\[0.2cm]
z'(t) = O(1) \quad \mathrm{as \ } t \downarrow 0^+ \quad , \quad
z(0^+) = z_0>0
\end{array}\right.
\end{equation}
Then
\begin{equation}\label{1/f}
\frac{1}{v} \in L^1(R^-)
\end{equation}
and for every $0<T<t<R$ such that $A\not\equiv 0$ in
$L^\infty([0,T])$
\begin{equation}\label{disugzero}
\int_T^t \big(\sqrt{A(s)}-\sqrt{\chi_R(s)}\big)ds \le -\frac{1}{2}
\Big( \log \int_0^TA(s)v(s)ds + \log \int_T^R \frac{ds}{v(s)}
\Big)
\end{equation}

\end{Teorema}
\noindent \textbf{Proof. } We set
\begin{equation} \label {uno:cinque}
y(t)=-\frac{v(t) z'(t)}{z(t)} \qquad \text{on } (0,R)
\end{equation}
Then $y \in \mathrm{Lip}_\mathrm{loc}([0,R))$; this follows since
$(vz')'= -Avz \in L^\infty_\mathrm{loc}([0,R))$, therefore $vz'$
is locally Lipschitz. Moreover, from \eqref{V} and \eqref{cauchy2}
we deduce that $y(0^+)=0$. Differentiating, we can argue that
$y(t)$ satisfies Riccati equation
\begin{equation} \label{uno:sei}
y'=A(t) v(t) + \frac{1}{v(t)} y^2 \qquad \text{a.e. on }  (0,R)
\end{equation}
Note that, since $A(t) \not \equiv 0$, $z$ is non constant and $y
\not\equiv 0$. Moreover, $y'(t) \ge 0$ almost everywhere on
$(0,R)$. From \eqref{A} and \eqref{uno:sei} it follows that, for
every $T
>0$ such that $A\not \equiv 0$ on $[0,T]$
\begin{equation}\label{uno:sette}
y(t) \ge y(T) \ge \int_0^T A(s)v(s)ds > 0\qquad \forall \ t \in
[T,R)
\end{equation}
From \eqref{uno:sei} and the elementary inequality $\epsilon a^2 +
\epsilon^{-1}b^2\ge 2 |a||b|$, $a,b \in\R$, $\epsilon >0$ we also
deduce $y' \ge 2 \sqrt{A(t)} |y(t)|$ and therefore
\begin{equation} \label{uno:otto}
y' \ge 2 \sqrt{A(t)} y \qquad \text{a.e. on }  [T,R)
\end{equation}
From \eqref{uno:sette} and \eqref{uno:otto} we infer
\begin{equation} \label{uno:nove}
\displaystyle y(t) \ge \Big(\int_0^T A(s)v(s)ds\Big) {\rm e}^{2
\int_{T}^t {\sqrt{A(s)} ds}} \qquad \text{on } [T,R)
\end{equation}
Moreover, from \eqref{uno:sei} and \eqref{A},
\begin{equation} \label{uno:undici}
\frac{y'}{y^2} \ge \frac{1}{v(t)} \qquad \text{a.e. on }  [T,R)
\end{equation}
Integrating on $[t,R-\varepsilon]$ for some small $\varepsilon>0$
we get
\begin{equation}\label{aaa}
\frac{1}{y(t)} \ge \frac{1}{y(R-\varepsilon)}
 + \int_t^{R-\varepsilon} \frac{ds}{v(s)} \ge \int_t^{R-\varepsilon}
\frac{ds}{v(s)}
\end{equation}
Letting $\varepsilon \rightarrow 0^+$ we obtain (\ref{1/f}), and
using (\ref{aaa}) into (\ref{uno:nove}) we reach the following
inequality:
\begin{equation}\label{primadis}
\int^t_{T} {\sqrt{A(s)} ds} \le -\frac{1}{2} \log \int_0^T
A(s)v(s)ds -\frac{1}{2} \log \int_t^{R} {\frac{ds}{v(s)}}
\end{equation}
Inequality (\ref{disugzero}) is simply a rewriting of
(\ref{primadis}): it is enough to point out that
\begin{equation}\label{intchi}
-\frac{1}{2} \log \int_t^{R} {\frac{ds}{v(s)}} = -\frac{1}{2} \log
\int_T^{R} {\frac{ds}{v(s)}} + \int_T^t \sqrt{\chi_R(s)}ds
\end{equation}
which follows integrating the definition of $\chi_R(t)$. $\blacktriangle$\\
\par
Although very simple, inequality (\ref{disugzero}) is deep. As we
have already stressed in the Introduction, the right hand side of
(\ref{disugzero}) is independent both of $t$ and of the behavior
of $A$ after $T$: if (\ref{disugzero}) is contradicted for some
$0<T<t<R$, the left hand side represents how much must $A(t)$
exceed the critical curve on the compact region $[T,t]$ in order
to have a first zero of $z(t)$, and it only depends on the
behavior of $A$ and $v$ before $T$ (the first addendum of the
right hand side), and on the growth of $v$ after $T$.\par
For geometrical purposes, from now on we will focus on the case
$R=+\infty$. However, the next corollaries can be restated on
$(0,R)$ replacing $+\infty$ with $R$ and $\chi(t)$ with
$\chi_R(t)$.
\begin{Remark}
\emph{Consider (\ref{intchi}) with $R=+\infty$:
$$
-\frac{1}{2} \log \int_t^{+\infty} {\frac{ds}{v(s)}} =
-\frac{1}{2} \log \int_T^{+\infty} {\frac{ds}{v(s)}} + \int_T^t
\sqrt{\chi(s)}ds \quad ,
$$
valid for $1/v \in L^1(+\infty)$. Letting $t\rightarrow +\infty$
we deduce that
\begin{equation}\label{noninl1}
\sqrt{\chi(t)} \not\in L^1(+\infty)
\end{equation}}
\end{Remark}
\begin{Corollario}[Existence of a first zero]\label{cor:zero}
In the assumptions of Theorem \ref{teo:uno} with $R=+\infty$,
suppose that either $1/v\not\in L^1(+\infty)$ or otherwise there
exist $0<T<t$ such that
\begin{equation}\label{contrad}
\int^{t}_{T} \big(\sqrt{A(s)}-\sqrt{\chi(s)}\big) ds >
-\frac{1}{2} \Big(\log \int_0^T A(s)v(s)ds +\log
\int_{T}^{+\infty} {\frac{ds}{v(s)}}\Big)
\end{equation}
Then, for every solution $z(t)\in
\mathrm{Lip}_\mathrm{loc}([0,+\infty))$ of \eqref{cauchy2}, there
exists $T_0 = T_0 (z)
>0$ such that $z(T_0) =0$. Moreover, the first
zero is attained on $(0,\overline{R}]$, where $\overline{R}>0$ is
the unique real number satisfying
\begin{equation}\label{prim}
\int^{t}_{T} {\sqrt{A(s)} ds} = -\frac{1}{2} \log \int_0^T
A(s)v(s)ds -\frac{1}{2} \log \int_{t}^{\overline{R}}
{\frac{ds}{v(s)}}
\end{equation}
\end{Corollario}
\noindent \textbf{Proof. } Observe that (\ref{contrad}) is
equivalent to say that (\ref{disugzero}) with $R=+\infty$ is false
for some $0<T<t$. Hence, the existence of a first zero on
$(0,+\infty)$ is immediate from Theorem \ref{teo:uno}.\par
As for the position of $T_0$, note first that (\ref{disugzero}) is
a rewriting of (\ref{primadis}). Suppose that $1/v \in
L^1(+\infty)$. We note that the $RHS$ of (\ref{primadis}) is
strictly decreasing as a function of $R\in (t,+\infty)$,
$\lim_{R\rightarrow t^-} RHS = +\infty$, and (\ref{primadis}) is
contradicted for $R=+\infty$ by assumption (\ref{contrad}).
Therefore, there exists a unique $\overline{R} \in (t,+\infty)$
such that (\ref{prim}) holds. Choosing $\varepsilon>0$ and
applying Theorem \ref{teo:uno} on the interval
$(0,\overline{R}+\varepsilon)$ we deduce the existence of a first
zero on $(0,\overline{R}+\varepsilon)$. Letting
$\varepsilon\rightarrow 0$ we reach the desired conclusion.\par
The case $1/v \not \in L^1(+\infty)$ is similar: we restrict the
considerations on a finite interval $[0,R]$, with $R>t$ small
enough that (\ref{primadis}) holds on $[0,R]$. Then, we enlarge
$R$ in such a way to reach the equality in (\ref{primadis}), and
we conclude as in the previous case. $\blacktriangle$\\
\begin{Corollario}[Oscillatory behavior] \label{corforte}
Fix $t_0 \in (0,+\infty)$. Suppose that \eqref{A}, \eqref{V} are
met on $[t_0,+\infty)$, with $1/v \in
L^\infty_\mathrm{loc}([t_0,+\infty))$, and let $z_0 \in
\mathbb{R}\backslash \{0\}$. Assume that either
\begin{equation}\label{casononint}
\frac{1}{v(t)} \not\in L^1(+\infty) \quad ,
\quad A(t)v(t) \not \in L^1(+\infty)
\end{equation}
or
\begin{equation}\label{condforte}
\frac{1}{v(t)} \in L^1(+\infty) \quad , \quad
\limsup_{t\rightarrow +\infty} \int_T^t
\big(\sqrt{A(s)}-\sqrt{\chi(s)}\big) ds = +\infty
\end{equation}
for some (hence any) $T>t_0$. Then, every solution $z(t)\in
\mathrm{Lip}_\mathrm{loc}([t_0,+\infty))$ of
\begin{equation}\label{cauchy2inf}
\left\{ \begin{array}{l} (v(t)z'(t))' + A(t)v(t) z(t) = 0 \qquad
\text{a.e. on } (t_0,+\infty)\\[0.2cm]
z(t_0) = z_0
\end{array}\right.
\end{equation}
is oscillatory.
\end{Corollario}
\noindent \textbf{Proof. } First, we claim that the two conditions
in \eqref{condforte} imply that $A(t)v(t) \not\in L^1(+\infty)$.
Indeed, from \eqref{noninl1} and the second condition of
\eqref{condforte} it follows that $\sqrt{A(t)} \not \in
L^1(+\infty)$, and from Cauchy-Schwarz inequality
$$
\Big(\int_T^t A(s)v(s)\Big)\Big(\int_T^t\frac{ds}{v(s)}\Big) \ge \Big(\int_T^t \sqrt{A(s)}ds\Big)^2
$$
letting $t\rightarrow +\infty$ we deduce the claim.\\
Suppose by contradiction that $z(t)$ has eventually constant sign.
Up to replacing $z$ with $-z$, we can assume $z(t)
> 0$ on $[\tau,+\infty)$, for some $\tau\ge t_0$. We define $y$ as in \eqref{uno:cinque}.
Then $y \in \mathrm{Lip}_\mathrm{loc}([\tau,+\infty))$ and
satisfies \eqref{uno:sei}, hence it is increasing. Integrating we
get
\begin{equation}\label{one}
y(t) \ge y(T) \ge y(\tau) + \int^T_\tau A(s)v(s)ds \qquad \forall
\ t>T>\tau
\end{equation}
By assumption, in both cases the non integrability of $A(t)v(t)$
ensures that there exists $T>\tau$ such that
$$
y(\tau) + \int_\tau^T A(s)v(s)ds > 0,
$$
therefore $y>0$ on $[T,+\infty)$. Now, we argue as in Theorem
\ref{teo:uno}. In particular, integrating \eqref{uno:undici} on
$[t,R_0]$ we get
\begin{equation}\label{two}
\frac{1}{y(t)} \ge \frac{1}{y(t)} - \frac{1}{y(R_0)} \ge
\int_t^{R_0} \frac{ds}{v(s)} \qquad \forall \ R_0>t>T
\end{equation}
so that $1/v \in L^1(+\infty)$, which contradicts
\eqref{casononint}. As for \eqref{condforte}, from $y' \ge
2y\sqrt{A}$ almost everywhere we deduce
\begin{equation}\label{three}
y(t) \ge y(T)\exp\left\{2\int_T^t \sqrt{A(s)} ds\right\} \qquad
\forall \ t>T.
\end{equation}
Combining \eqref{one}, \eqref{two}, \eqref{three} and using the
definition of $\chi(t)$ we obtain the following inequality:
$$
\int^t_{T} \big(\sqrt{A(s)}-\sqrt{\chi(s)}\big) ds \le -\frac{1}{2} \log \Big(y(\tau)+\int_\tau^T
A(s)v(s)ds\Big) -\frac{1}{2} \log \int_T^{+\infty} {\frac{ds}{v(s)}}
$$
Letting $t\rightarrow +\infty$ along a sequence realizing
(\ref{condforte}) we reach the desired
contradiction. $\blacktriangle$\\

\par
Here are some stronger conditions which imply oscillation, and
that will be used in the sequel.
\begin{Proposizione}\label{condsuff}
In the assumptions \eqref{A}, \eqref{V} on the interval
$[t_0,+\infty)$, equation \eqref{cauchy2inf} is oscillatory in the
following cases:
\begin{itemize}
\item[-] $1/v \in L^1(+\infty)$ and one of the following
    conditions is satisfied for some $T>t_0$:
$$
\begin{array}{cl}
(i) & \displaystyle A(t) \ge \chi(t) \quad \text{a.e. on } [T,+\infty) \quad
\text{and} \quad \sqrt{A(s)}-\sqrt{\chi(s)} \not\in
L^1(+\infty) \\[0.2cm]
(ii) & \displaystyle \limsup_{t\rightarrow +\infty} \frac{\int_T^t
\sqrt{A(s)}ds}{\int_T^t \sqrt{\chi(s)}ds} >1 \\[0.5cm]
(iii) & \displaystyle \liminf_{t\rightarrow +\infty}
\frac{\sqrt{A(t)}}{\sqrt{\chi(t)}} > 1 \\[0.5cm]
(iv) & \displaystyle \limsup_{t\rightarrow +\infty} \frac{\int_T^t
\sqrt{A(s)}ds}{ -\frac{1}{2} \log \int_t^{+\infty}
\frac{ds}{v(s)}} > 1
\end{array}
$$
\item[-] $v(t) \not\in L^1(+\infty)$, $v(t) \le f(t)$ a.e. for
some continuous function $f(t)$ such that $1/f \in L^1(+\infty)$
and
$$
\begin{array}{cl}
(v) & A \ \mathit{is \ positive, \ increasing \ and \ }
\displaystyle \sqrt{A(t_n)} > \inf_{t > t_n} \Big\{ -\frac{1}{2}
\frac{\log \int_t^{+\infty}
\frac{ds}{f(s)}}{t-t_n} \Big\} \\[0.3cm]
& \mathit{for \ some \ increasing \ sequence } \  \{t_n\} \uparrow
+\infty
\end{array}
$$
\end{itemize}
\end{Proposizione}
\noindent \textbf{Proof. } Implications $(i),(ii),(iii)$ are
immediate from (\ref{noninl1}). To obtain $(iv)$ we also use
equality (\ref{intchi}) with $R=+\infty$. Regarding $(v)$, we
proceed, by contradiction, as in Corollary \ref{corforte},
restricting the problem on $[\tau,+\infty)$, $\tau>t_0$. Since
$A(t)$ is increasing, it is bounded from below away from zero on
$[\tau,+\infty)$. Therefore, since $v(t) \not\in L^1(+\infty)$ we
can choose $T>\tau$ such that
$$
y(\tau) + \int_\tau^T A(s) v(s)ds \ge 1.
$$
Using the monotonicity of $A$ and $v \le f$, (\ref{primadis})
becomes
$$
\sqrt{A(T)}(t-T) \le \int_T^t \sqrt{A(s)}ds \le -\frac{1}{2}
\log\int_t^{+\infty} \frac{ds}{v(s)} \le -\frac{1}{2}
\log\int_t^{+\infty} \frac{ds}{f(s)}
$$
for every $T<t$; $(v)$ contradicts this last chain of
inequalities. $\blacktriangle$\\
\par
Corollary \ref{corforte} is related to the classical Hille-Nehari
oscillation theorem (see \cite{S}). However, in order to apply
this latter to ensure that a solution $z(t)$ of (\ref{cauchy2inf})
is oscillatory, one needs to perform a change of variables which
requires $1/v \in L^1(+\infty)$. Therefore, Hille-Nehari criterion
is not straightforwardly applicable when $1/v \not\in
L^1(+\infty)$. Moreover, in case $1/v \in L^1(+\infty)$, in order
to have oscillatory solutions the criterion requires that
\begin{equation}\label{HN}
\liminf_{t\rightarrow +\infty} \sqrt{A(t)} v(t) \int_t^{+\infty}
\frac{ds}{v(s)} > \frac{1}{2}
\end{equation}
which is exactly request $(iii)$ of Proposition \ref{condsuff},
using definition (\ref{chi}) of $\chi(t)$. It is worth to point
out that (\ref{condforte}) implies oscillations even in some cases
when the ``liminf" in (\ref{HN}) is equal to $1/2$, an
unpredictable case in Hille-Nehari theorem.
\section{Why is the critical curve really critical?}
In this section we show that Corollary \ref{corforte} is sharp.
This will be done by studying the relationship between $\chi(t)$
and
the two critical functions introduced in \cite{BR}.\\
\par
Consider the ``Euclidean" problem
\begin{equation}\label{cauchyeuclidea}
\left\{ \begin{array}{l} (t^{m-1}z'(t))' + A(t)t^{m-1} z(t) = 0
\qquad \mathrm{on \ }
(0,+\infty) \\[0.1cm]
z'(0^+) = 0 \quad , \quad z(0) = z_0>0 \quad , \quad m\ge 3
\end{array}\right.
\end{equation}
In this case, from $v(t)= t^{m-1}$ it is immediate to see that
\begin{equation}\label{chieu}
\chi(t) = \frac{(m-2)^2}{4}\frac{1}{t^2}.
\end{equation}
Suppose that $0\le A(t) \in C^\infty([0,+\infty))$ is such that,
for some $\varepsilon>0$,
\begin{equation}\label{euc}
A(t) \left\{\begin{array}{ll} \le \displaystyle
\frac{(m-2)^2}{4}\frac{1}{t^2} &
\mathrm{on \ } [0,\varepsilon) \\[0.3cm]
\displaystyle = \frac{(m-2)^2}{4}\frac{1}{t^2} & \mathrm{on \ }
[\varepsilon,+\infty)
\end{array}\right.
\end{equation}
Then, problem (\ref{cauchyeuclidea}) admits a positive solution
$0<z(t) \in C^1([0,+\infty))$ satisfying, by Proposition $4.1$ of
\cite{BR},
$$
C^{-1} t^{-\frac{m-2}{2}} \log t \le z(t) \le C t^{-\frac{m-2}{2}}
\log t
$$
for some positive constant $C$ and $t\gg 1$. Suppose now that
$A(t) = H^2/t^2$ on $[\varepsilon,+\infty)$. By Proposition
\ref{comparison} in Appendix, there exists a positive solution for
every $H \le \frac{m-2}{2}$, while in case $H>\frac{m-2}{2}$ the
limit in item $(iii)$ of (\ref{condsuff}) is
\begin{equation}\label{limite}
\lim_{t \rightarrow +\infty } \frac{\displaystyle
\sqrt{A(t)}}{\displaystyle \sqrt{\chi(t)}} = \frac{2H}{m-2} > 1,
\end{equation}
and by Corollary \ref{corforte} every solution $z(t)$ is
oscillatory. Therefore, in the Euclidean case we recognize
(\ref{chieu}) as the correct critical curve for the behavior of $z(t)$.\\
\par
The hyperbolic case is less immediate. However, fix $B>0$ and
consider
\begin{equation}\label{cauchyiperbolica}
\left\{ \begin{array}{l} (\sinh^{m-1}(Bt)z'(t))' +
A(t)\sinh^{m-1}(Bt) z(t) = 0 \qquad \mathrm{on \ }
(0,+\infty) \\[0.1cm]
z'(0^+) = 0 \quad , \quad z(0) = z_0>0 \quad , \quad m\ge 2
\end{array}\right.
\end{equation}
In this case $v(t)= \sinh^{m-1}(Bt)$ and the expression of
$\chi(t)$ is more complicated. Nevertheless, using De l'Hopital
theorem, we see that, as $t\rightarrow +\infty$,
$$
\displaystyle \chi(t) = \displaystyle
\Big[\frac{1}{2\sinh^{m-1}(Bt)\int_t^{+\infty}
\sinh^{1-m}(Bs)ds}\Big]^2  \sim \displaystyle \frac{(m-1)^2B^2}{4}
\coth(Bt)
$$


Suppose now that $0 \le A(t) \in C^\infty([0,+\infty))$ is such
that, for some $\varepsilon>0$,
\begin{equation}\label{hyp}
A(t) \left\{\begin{array}{ll} \le \displaystyle
\frac{(m-1)^2B^2}{4}\coth (Bt) &
\mathrm{on \ } [0,\varepsilon) \\[0.3cm]
\displaystyle = \frac{(m-1)^2B^2}{4}\coth (Bt) & \mathrm{on \ }
[\varepsilon,+\infty)
\end{array}\right.
\end{equation}
Then, (\ref{cauchyiperbolica}) has a positive solution $z\in
C^1([0,+\infty))$ satisfying
$$
C^{-1} t \mathrm{e}^{-\frac{m-1}{2}Bt} \le z(t) \le C t
\mathrm{e}^{-\frac{m-1}{2}Bt}
$$
for some appropriate constant $C>0$ and $t\gg 1$.
\par
In case $A(t)= H^2B^2\coth (Bt)$ on $[\varepsilon,+\infty)$, again
using Proposition \ref{comparison} we deduce that, for every $H\le
\frac{m-1}{2}$, there exists a positive solution of
(\ref{cauchyiperbolica}). On the contrary, if $H>\frac{m-1}{2}$
the limit in item $(iii)$ of Proposition \ref{condsuff} is
strictly greater than $1$, hence every solution is oscillatory.
The characteristic curve $\chi(t)$ is ``asymptotically sharp" even
in the hyperbolic case, and numerical evidences show it agrees
sharply with the curve $\frac{(m-1)^2B^2}{4} \coth(Bt)$ outside
$t=0$.
\section{Oscillation estimates: the key
result}\label{oscillestimates}
So far, we have only ensured an oscillatory behavior of solutions
of (\ref{cauchy2inf}) in case $A(t)$ is, for example, asymptotic
to the critical curve and $\sqrt{A(t)}-\sqrt{\chi(t)}$ is
eventually positive and non integrable at infinity. Under these
assumptions, we cannot expect the oscillations to be automatically
thick, since we have proved that $\chi(t)$ is sharp as a border
line function. Nevertheless, suppose that
$$
\frac{A(t)}{\chi(t)} \ge c>1 \qquad \mathrm{for} \ t\gg 1.
$$
In this case, one may expect that the somewhat ``uniform" mass of
$A(t)$ exceeding from $\chi(t)$ can control the distance between
zeros from above. The key result, Theorem \ref{teo:princ}, goes in
this direction: given two consecutive zeros $T_1(\tau)<T_2(\tau)$
of $z(t)$ after $\tau$ it states that
$$
T_2(\tau) - T_1(\tau) = O(\tau) \qquad \mathrm{as \ }
\tau\rightarrow +\infty.
$$
Moreover, in case
\begin{equation}\label{protot}
v(t) \le f(t) = \Lambda \exp\{at^\alpha
\log^\beta t\} \qquad \Lambda,a,\alpha >0 \ , \ \beta \ge 0
\end{equation}
we will be able to estimate the quantity
$$
\limsup_{\tau \rightarrow +\infty} \frac{T_2(\tau)}{\tau}.
$$
Theorem \ref{teo:princ} exploits upper bounds for the function
$v(t)$ in terms of some function $f(t)$, instead of dealing with
$v(t)$ itself. The necessity of working with such an upper bound
needs some preliminary comment.
\par
Although the critical function $\chi(t)$ is suitable to describe
the oscillatory behavior of (\ref{cauchy2inf}), due to its
integral expression in $v(t)$ it is in general not easy to handle.
Moreover, $v(t)$ itself can behave very badly since, in our
geometric applications, it represents the volume growth of
geodesic spheres; indeed, in many situations, such as volume
comparison results, one deals only with upper bounds of the volume
growth in terms of some known function $f(r)$ which possesses some
further regularity property (for example, as we will suppose in
the sequel, monotonicity and differentiability). Hence, it would
be useful to look for slightly precise but more manageable
critical functions depending on $f(t)$ instead of $v(t)$. The most
natural way is to define
\begin{equation}\label{chif}
\chi_f(t) = \Big[\frac{1}{2f(t)\int_t^{+\infty}
\frac{ds}{f(s)}}\Big]^2 = \Big[\Big( -\frac{1}{2}
\log\int_t^{+\infty} \frac{ds}{f(s)} \Big)'\Big]^2.
\end{equation}
Since the integral of $\sqrt{\chi_f(t)}$ is greater than the
integral of $\sqrt{\chi(t)}$ on every compact interval, we have
that $\chi_f(t) \ge \chi(t)$ on $(0,+\infty)$, and obviously
$\chi_f\equiv \chi$ in case $v\equiv f$. It is not hard to see
that, if we substitute $\chi(t)$ with $\chi_f(t)$ and $v(t)$ with
$f(t)$ (with the exception of the terms involving integrals of
$A(t)v(t)$), all the conclusions of the theorems of Section
\ref{existencefirstzero} are still true.\par
Unfortunately, despite the further properties of $f$, even this
critical function is too difficult to handle in many instances.
Hence, we choose the simpler critical function
\begin{equation}\label{chifsemplice}
\widetilde{\chi}_f(t) = \Big[\frac{f'(t)}{2f(t)}\Big]^2.
\end{equation}
Since (\ref{protot}) represents the prototype of most volume
growth bounds, it is important to stress the relationship between
$\chi_f(t)$ and $\widetilde{\chi}_f(t)$ in case $f(t)=\Lambda
\exp\{at^\alpha \log^\beta t\}$. Using De l'Hopital theorem we
have
\begin{equation}\label{asinttt}
\lim_{t\rightarrow +\infty}
\frac{\sqrt{\widetilde{\chi}_f(t)}}{\sqrt{\chi_f(t)}} =
\lim_{t\rightarrow +\infty} \frac{f'(t)^2}{f(t) f''(t)} = 1 \qquad
\mathrm{since \ } \alpha>0.
\end{equation}
Therefore, with this choice of $f$ the modified critical function
$\widetilde{\chi}_f(t)$ is asymptotic to the critical function
$\chi_f(t)$. This justifies the use of $\widetilde{\chi}_f(t)$ as
a border line ``at infinity" for $A(t)$.
\par
Throughout this section we shall require the validity of the
following properties on $[t_0,+\infty)$, for some $t_0>0$.
\begin{align}
&0 \le v(t) \in L^\infty_\mathrm{loc}([t_0,+\infty)) , \quad
\frac{1}{v(t)} \in L^\infty_\mathrm{loc}([t_0,+\infty)) , \quad
\frac{1}{v(t)} \in
L^1(+\infty) & & \tag{V2}\label{v1s} \\[0.1cm]
&f \in C^1 ([t_0,+\infty)) \quad , \quad f(t_0) >0 \quad  & & \tag{F1}
\label{f1p}
\\[0.1cm]
&f \ {\rm is \ non \ decreasing \ on  } \ [t_0,+\infty) & & \tag{F2} \label{f6}\\[0.1cm]
& v(t) \le f(t) \quad \text{a.e. on }  [t_0,+\infty) & & \tag{F3}
\label{f7}
\\[0.1cm]
& \forall \ t \ge t_0 \quad \frac{f'(t)}{f(t)} \ge \frac{1}{Dt^\mu} \quad {\rm for \ some} \ D>0, \ \mu <1 & & \tag{F4} \label{F5}\\[0.1cm]
& A \in L^\infty_\mathrm{loc}([t_0,+\infty)), \quad A(t)\ge 0 \quad \text{a.e. on } [t_0,+\infty) & & \tag{A2}\label{a1s} \\[0.1cm]
& \limsup_{t\rightarrow +\infty} \int_{t_0}^t \big(\sqrt{A(s)} -
\sqrt{\chi(s)}\big)ds = +\infty & & \tag{A3} \label{A5} \\[0.1cm]
&\exists \ c>0 \ {\rm such \ that \ } \quad \sqrt{A(t)} \ge c
\sqrt{\widetilde{\chi}_f(t)} = \frac{c}{2}\frac{f'(t)}{f(t)}
\qquad \text{a.e. on } [t_0,+\infty) & & \tag{A4} \label{A4}
\end{align}\\
Next, we introduce two classes of functions: for $f \in
C^0([t_0,+\infty))$, $f>0$ on $[t_0,+\infty)$, $h, k$ piecewise
$C^0$ and non-negative on $[t_0,+\infty)$, $c>0$ we set
\begin{equation} \label{tre:uno}
\begin{array}{ll}
\displaystyle {\cal A}(f,h,c) = & \Big\{ g  :  [t_0,+\infty) \to [0,+\infty) \ {\rm \ piecewise \ } C^0 {\rm \ such \ that}\\
&  \left. \displaystyle \limsup_{t \to +\infty} \left( \sup_{\xi
\in (0,1)} \frac{(1-\xi) g(t)
f(t+g(t)+h(t))^c}{f(t+(1-\xi)g(t)+h(t))^{c+1}} \right) <+\infty
\right\}
\end{array}\\[0.3cm]
\end{equation}
\begin{equation}\label{tre:due}
\begin{array}{ll}
\displaystyle {\cal B}(f,k,c) = & \Big\{ g :  [t_0,+\infty) \to [0,+\infty) \ {\rm \ piecewise \ } C^0  {\rm \ such \ that} \\
& \left. \displaystyle \limsup_{t \to +\infty} \left( \sup_{\xi
\in (0,1)} \frac{\xi g(t) f(t+(1-\xi)g(t)+k(t))^c}{f(t+g(t)+k(t))
\cdot f(t+k(t))^c}\right) <+\infty  \right\}
\end{array}\\[0.4cm]
\end{equation}
\textbf{Definition:} \emph{We  shall say that $f$
\emph{\textbf{satisfies property $(P)$ for some $c>0$}} if
whenever $$ h(t), \ k(t) = O(t) \quad \mathit{as} \ t\to +\infty
\quad , \quad g \in {\cal A}(f,h,c) \cup
{\cal B} (f,k,c) $$ implies $g(t)= O(t)$  as $t \to +\infty$.}\\
%
\par
An example of $f$ satisfying property $(P)$ that we shall use in
the sequel is the following. Let
\begin{equation}\label{prototipomain}
f(t)=\exp\{a t^\alpha \log
^{\beta} t\} \qquad a>0, \ \alpha
> 0, \ \beta \ge 0 \qquad {\rm for} \ t \ge t_0.
\end{equation}
Then \emph{$f$ satisfies property $(P)$ for every $c>1$}. Indeed,
let $h$ and $k$ be non negative and such that $h(t), \ k(t)=O(t)$
as $t \to +\infty$ and let  $g \in  {\cal A}(f,h,c)$. Assume, by
contradiction, the existence of a sequence  $\{t_n\} \to +\infty$
with the property
\begin{equation} \label{tre:cinque}
\frac{g(t_n)}{t_n} \to +\infty \qquad {\rm as} \ n \to +\infty
\end{equation}
Without loss of generality we suppose $g(t_n) >1$ $\forall \ n$
and we define $\xi_n=1 -\frac{1}{g(t_n)}$. Then
\begin{align}
& \displaystyle
\frac{(1-\xi_n)g(t_n)f(t_n+g(t_n)+h(t_n))^c}{f(t_n+(1-\xi_n)g(t_n)+h(t_n))^{c+1}}
= \frac{f(t_n+g(t_n)+h(t_n))^c}{f(t_n+1+h(t_n))^{c+1}} = & &
\label{casoa}
\\[0.2cm]
& \displaystyle = \exp\Big\{ac(t_n+g(t_n)+h(t_n))^\alpha
\log^\beta
(t_n+g(t_n)+h(t_n)) \ + & & \nonumber \\
& - a(c+1)(t_n+1+h(t_n))^\alpha \log^\beta
(t_n+1+h(t_n))\Big\} = & & \nonumber \\
& = \exp\Big\{acg(t_n)^\alpha \log^\beta (t_n+g(t_n)+h(t_n))\cdot
& & \nonumber \\
&  \displaystyle \Big[\Big(1 + \frac{t_n}{g(t_n)} +
\frac{h(t_n)}{g(t_n)}\Big)^\alpha + & & \label{primo}
\\
& \displaystyle  - \frac{(c+1)t_n^\alpha}{cg(t_n)^\alpha} \Big(1 +
\frac{1}{t_n} + \frac{h(t_n)}{t_n}\Big)^\alpha \frac{\log^\beta
(t_n+1+h(t_n))}{\log^\beta (t_n+g(t_n)+h(t_n))} \Big]\Big\}  & &
\label{secondo}
\end{align}
Note that expression \eqref{primo} tends to $1$ as $n\to+\infty$,
while expression \eqref{secondo} goes to $0$. Their difference is
thus eventually positive, so \eqref{casoa} goes to $+\infty$, but
this contradicts the fact that $g \in {\cal A}(f,h,c)$. Observe
that here any $c>0$ would work. Let now $g \in {\cal B}(f,k,c)$
and reason again by contradiction. Let $\{t_n\}$ be as above. Then
\begin{align}
& \displaystyle \xi g(t_n)\frac{f(t_n+ (1-\xi)g(t_n) +
k(t_n))^c}{f(t_n+g(t_n) +k(t_n))\cdot f(t_n+k(t_n))^c} = & &
\label{casob}\\[0.2cm]
& = \xi g(t_n) \exp \Big\{ac(1-\xi)^\alpha g(t_n)^\alpha \Big(1 +
\frac{1}{1-\xi}\Big(\frac{t_n}{g(t_n)}+\frac{k(t_n)}{g(t_n)}\Big)\Big)^\alpha
\cdot & & \nonumber \\
& \log^\beta(t_n+(1-\xi)g(t_n) +k(t_n)) - ag(t_n)^\alpha \Big(
1+\frac{t_n}{g(t_n)}+\frac{k(t_n)}{g(t_n)}\Big)^\alpha \cdot & &
\nonumber
\\
& \log^\beta (t_n+g(t_n)+k(t_n)) - ac t_n^\alpha \Big(
1+\frac{k(t_n)}{t_n}\Big) ^\alpha \log^\beta (t_n+k(t_n)) \Big\} =
& & \nonumber \\
& \ge \xi g(t_n) \exp \Big\{a g(t_n)^\alpha
\log^\beta(t_n+(1-\xi)g(t_n) +k(t_n)) \cdot & & \nonumber \\
&  \Big[\Big(c(1-\xi)^\alpha -
\frac{\log^\beta(t_n+g(t_n)+k(t_n))}{\log^\beta(t_n+(1-\xi)g(t_n)+k(t_n))}\Big)
\Big(
1+\frac{t_n}{g(t_n)}+\frac{k(t_n)}{g(t_n)}\Big)^\alpha \ + & & \label{terzo} \\
&  \displaystyle - c
\frac{t_n^\alpha}{g(t_n)^\alpha}\Big(1+\frac{k(t_n)}{t_n}\Big)^\alpha
\frac{\log^\beta(t_n+k(t_n))}{\log^\beta(t_n + (1-\xi)g(t_n) +
k(t_n))} \Big]\Big\} \label{quarto}
\end{align}
Since expression \eqref{quarto} goes to $0$ as $n\to+\infty$, we
can choose $n$ such that it is eventually less than $\epsilon$,
for some fixed $\epsilon>0$. Moreover, since $\forall \ \xi \in
(0,1)$
$$
\displaystyle
\frac{\log^\beta(t_n+g(t_n)+k(t_n))}{\log^\beta(t_n+(1-\xi)g(t_n)+k(t_n))}
\to 1 \quad {\rm as \ } n\to +\infty
$$
and using now $c>1$, we can choose a suitable $\xi$ such that
expression \eqref{terzo} is eventually strictly positive and
greater than $2\epsilon$, if we choose $\epsilon$ sufficiently
small. Now letting $n \to +\infty$ we have that \eqref{casob} goes
to infinity, which implies $ g \not\in {\cal B}(f,k,c)$, a
contradiction. Note that assumption $\alpha>0$ is necessary: it is
not hard to see that, if $f(t)$ has polynomial growth, then $f$
does not satisfy property $(P)$ for any $c>0$. On the contrary,
proceeding in a way similar to that outlined above one verifies,
for instance, that also the function
$$
\Lambda \exp\{a\mathrm{e}^{bt}\} \quad , \quad \Lambda,a,b>0
$$
satisfies property $(P)$ for every $c>1$. Assuming $f(t)$ of this
type, one can prove analogous estimates as those in
(\ref{limsupa}) and (\ref{tesiuno}).\\


\par
Going back to (\ref{prototipomain}), we observe that \eqref{f1p},
\eqref{f6} and \eqref{F5} are satisfied. We also observe that the
validity of \eqref{v1s}, \eqref{a1s} and \eqref{A5} enables us to
apply Corollary \ref{corforte} to conclude that equation
(\ref{cauchy2inf}) is oscillatory on $[t_0,+\infty)$, and that, by
Proposition \ref{isolati} of the Appendix, the zeros of $z(t)$ are
isolated.
\par
 Now, we are ready to prove our main technical result.
\begin{Teorema}\label{teo:princ} Assume the validity of
\eqref{v1s}, \eqref{f1p}, \eqref{f6}, \eqref{f7}, \eqref{F5},
\eqref{a1s}, \eqref{A5}, \eqref{A4} and that $f$ satisfies
property $(P)$ for the parameter $c>0$ required in \eqref{A4}. Let
$z \not\equiv 0$ be a locally Lipschitz solution of
\eqref{cauchy2inf} on $[t_0,+\infty)$. Let $\tau \in [T,+\infty)$,
where $T$ is defined in Corollary \ref{corforte}, and let
$T_1(\tau)$, $T_2 (\tau)$ be the first two consecutive zeros of
$z(t)$ on $[\tau,+\infty)$. Then
\begin{equation} \label{tre:sei}
T_2(\tau) -\tau= O(\tau) \qquad {\rm as} \ \tau \to  +\infty.
\end{equation}
Moreover, in case $f(t)= \Lambda \exp[at^\alpha \log^\beta t]$ we
have the estimate
\begin{equation}\label{limsupa}
\limsup_{\tau \to+\infty} \frac{T_2(\tau)}{\tau} \le
\Big(\frac{c+1}{c-1}\Big)^{\frac{2}{\alpha}}.
\end{equation}


\end{Teorema}
\noindent \textbf{Proof. } As we have observed, $z(t)$ is
oscillatory. Having fixed $\tau \in [T,+\infty)$, let
$$
U=[\tau,T_2(\tau))\backslash \{ T_1(\tau) \}
$$
and on $U$ consider the locally Lipschitz function
$$
y(t)=- \frac{v(t) z' (t)}{z(t)}
$$
solution of
\begin{equation}\label{tre:sette}
y'(t) = A(t) v(t) + \frac{1}{v(t)} y^2(t) \qquad \text{a.e. on }
[t_0,+\infty)
\end{equation}
Because of \eqref{a1s} and  \eqref{v1s}, \eqref{tre:sette} shows
that $y$ is non decreasing on $U$. Indeed, from \eqref{A4},
\eqref{F5}, \eqref{v1s} we can argue that $y$ is strictly
increasing on $U$. Since $z\not\equiv 0$, proceeding analogously
to Proposition \ref{isolati} in Appendix we deduce that
\begin{equation}\label{tre:otto}
y(T_1(\tau)^+)=-\infty, \quad y(T_1(\tau)^-)=+\infty, \quad
y(T_2(\tau)^-)=+\infty
\end{equation}
Note that it could be $U=(T_1 (\tau), T_2 (\tau))$: this is
exactly the case when $T_1(\tau)=\tau$. \par Due to the fact that
$y$ is non decreasing, $U$ can be decomposed as a disjoint union
of intervals of the types
\begin{equation}\label{livelli}
\begin{array}{ll}
I_1 \subseteq \{ x \in U \ : \ y(x) \in [-1,1]\}  & \quad {\rm
interval \ of \ type \ 1}  \\[0.1cm]
I_2 \subseteq \{ x \in U \ : \ y(x) >1\} & \quad {\rm interval \ of \ type \ 2}
\\[0.1cm]
I_3 \subseteq \{ x \in U \ : \ y(x) <-1\} & \quad {\rm interval \
of \ type \ 3}
\end{array}
\end{equation}
To fix ideas we consider the case $y(\tau) <-1$ pictured below.
\begin{figure}[ht]
       \includegraphics[width=12cm]{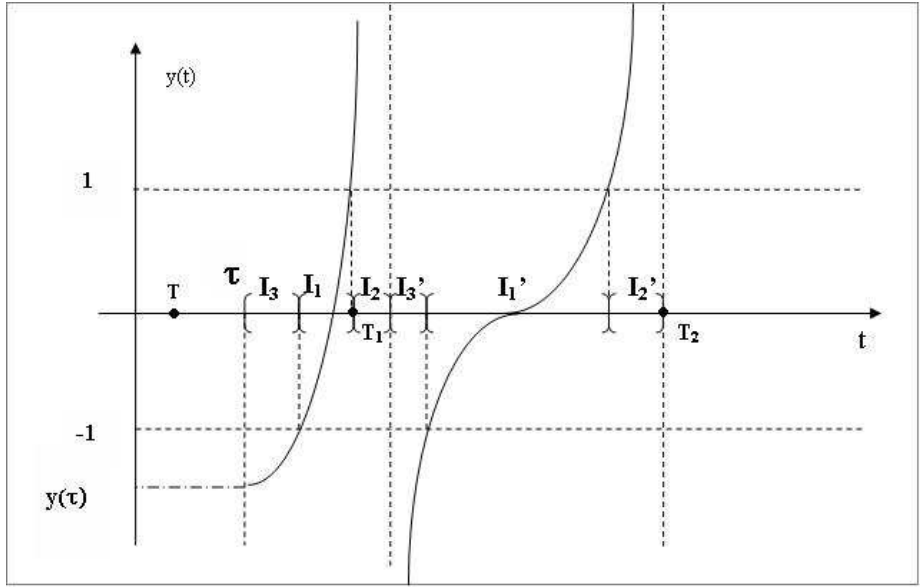}
       \caption{Riccati solution}
\end{figure}

In this case we have
$$
U=I_3 \cup I_1 \cup I_2 \cup I_3' \cup I_1' \cup I_2'
$$
where: \par
 $I_1$ is the first interval of type $1$, after $\tau$ and before $T_1 (\tau)$;\par
 $I_2$ is the first interval of type $2$, after $\tau$ and before $T_1(\tau)$;\par
 $I_3$ is the first interval of type $3$, after $\tau$ and before $T_1(\tau)$;\par
 $I_1'$ is the first interval of type $1$, after $T_1(\tau)$ and before $T_2 (\tau)$; \par
 $I_2'$ is the first interval of type $2$, after $T_1 (\tau)$ and before $T_2 (\tau)$;\par
 $I_3'$ is the first interval of type $3$, after $T_1 (\tau)$ and before $T_2(\tau)$.
 \\
 \\
  We study this situation which is ``the worst" it could happen. The remaining cases can be dealt with similarly and we shall skip proofs. \par
  For $i =\{1,2,3\}$ we set $|I_i|=g_i (\tau)$ and $|I_i'|=g_i '(\tau)$.
  We are going to prove that, in the above hypotheses, each $g_i(\tau), \ g_i'(\tau)$ is $O(\tau)$ as $\tau \to
  +\infty$.\\
  \par

We consider at first an open interval $J$ of type $3$ so that $J$
could be either $I_3$ or $I_3'$. Set $P(\tau) <Q(\tau)$ to denote
its end points; thus $g_3(\tau)=|J|(\tau)=Q(\tau)-P(\tau)$ and
$g_3(\tau)$ is clearly piecewise $C^0([T,+\infty))$. We have
$y(Q)=-1$ and $y(P) \le -1$ if $y$ is defined in $P$, otherwise
$y(P^+) = -\infty$. As in Theorem \ref{teo:uno}, \eqref{tre:sette}
yields
$$
y' \ge 2 \sqrt{A(t)} |y| = 2 \sqrt{A(t)} (-y) \qquad \text{a.e. on
} J
$$
Fix $t \in (P,Q]$ and integrate on $[t,Q]$. Recalling that $y(s)
\le y(Q)=-1$ $\forall \ s \in (P,Q]$ we have
\begin{equation}\label{tre:nove}
y(t) \le - \exp\Big\{2\int_t^Q {\sqrt{A(s)}ds}\Big\} \qquad \forall \ t\in
(P,Q]
\end{equation}
Since $y'/y^2 \ge 1/v$ almost everywhere, integrating on
$[P+\varepsilon,t]$ for some small $\varepsilon>0$ we obtain
\begin{equation}\label{tre:a}
\frac{1}{y(P+\varepsilon)}-\frac{1}{y(t)} \ge
\int^{t}_{P+\varepsilon} {\frac{ds}{f(s)}}
\end{equation}
Letting $\varepsilon \to 0^+$ we get
\begin{equation}\label{tre:dodici}
-\frac{1}{y(t)} \ge -\frac{1}{y(P^+)} + \int_P^{t} \frac{ds}{f(s)}
\ge \int_P^{t} \frac{ds}{f(s)}
\end{equation}
which is valid $\forall$ $t \in (P,Q]$. Now, from \eqref{tre:nove}
and because of \eqref{A4}
$$
2 \int_{t}^Q \sqrt{A(s)} ds \ge c \int_{t}^Q \frac{f'(s)}{f(s)} ds
= \log \left(\frac{f(Q)}{f(t)} \right)^c
$$
and therefore, from \eqref{tre:nove},
$$
-\frac{1}{y(t)} \le \left( \frac{f(t)}{f(Q)} \right)^c
$$
Substituting into \eqref{tre:dodici} and using \eqref{f6} we
obtain
\begin{equation}\label{tre:quattordici}
1 \ge \left( \frac{f(Q)}{f(t)} \right)^c \int_P^{t}
\frac{ds}{f(s)} \ge (t-P) \frac{f(Q)^c}{f(t)^{c+1}} \quad \forall
\ t \in (P,Q)
\end{equation}
Suppose now that $J=I_3$, so that $P(\tau)=\tau$ and
$Q(\tau)=\tau+g_3(\tau)$ . Since $t \in (P,Q)$, there exists $\xi
\in (0,1)$ such that
$$
t=\tau + (1-\xi) g_3(\tau) \quad , \quad t-P=(1-\xi)g_3(\tau)
$$
and since $t$ was arbitrary, from \eqref{tre:quattordici} we
obtain
\begin{equation}\label{gitre}
\sup_{\xi \in (0,1)} \frac{(1-\xi) g_3(\tau)
f(\tau+g_3(\tau))^c}{f(\tau+(1-\xi) g_3(\tau))^{c+1}} \le 1
\end{equation}
in this case it follows that $g_3 \in {\cal A}(f,0,c)$ and then
$g_3(\tau)=O(\tau)$ as $\tau \to +\infty$. \\
We will deal with the case $J=I_3'$ later.\\
\par
Next, we consider an interval $J$ of type $1$. Set
$P(\tau)<Q(\tau)$ to denote its end points; thus
$g_1(\tau)=|J|(\tau)=Q(\tau)-P(\tau)$ and $g_1(\tau)$ is piecewise
$C^0([T;+\infty))$. In this case $y(P)=-1$,  $y(Q)=1$ and $|y| \le
1$ on $J$. We integrate Riccati equation \eqref{tre:sette} on
$[P,Q]$ to obtain
$$
2=\int_P^Q {y'(s) ds}=\int_P^Q {A(s)v(s)ds}+\int_P^Q
\frac{y^2(s)}{v(s)} ds \ge \int_P^Q A(s) v(s) ds
$$
Next, without loss of generality we can suppose to have chosen $T$
sufficiently large that \eqref{v1s}, in particular  $1/v \in
L^1(+\infty)$, implies
$$
\int_T^{+\infty} \frac{ds}{v(s)} \le 1
$$
so that
$$
\int_P^Q \frac{ds}{v(s)} \le 1
$$
From the above inequality, using \eqref{A4} and the generalized
mean value theorem it follows that, for some $T_0 \in [P,Q]$,
$$
\begin{array}{lcl}
2 & \ge & \displaystyle \int_P^Q {A(s) v(s) ds} \int_P^Q
{\frac{ds}{v(s)}} \ge \int_P^Q
\frac{c^2}{4}\left(\frac{f'(t)}{f(t)}\right)^2 v(s) ds \int_P^Q
\frac{ds}{v(s)} \\[0.5cm]
& = & \displaystyle
\frac{c^2}{4}\left(\frac{f'(T_0)}{f(T_0)}\right)^2 \int_P^Q {v(s)
ds} \int_P^Q {\frac{ds}{v(s)}}
\end{array}
$$
On the other hand, from H\"older inequality
$$
(Q-P)^2 \le \int_P^Q {v(s) ds} \int_P^Q {\frac{ds}{v(s)}}
$$
so that
$$
\sqrt{2} \ge \frac{c}{2} \left( \frac{f'(T_0)}{f(T_0)} \right)
(Q-P)
$$
or, in other words, using \eqref{f1p}, \eqref{f6} and observing
that \eqref{F5} implies that $f'$ is eventually positive,
\begin{equation} \label {tre:venti}
\frac{2\sqrt{2}}{c} \frac{f(T_0)}{f'(T_0)} \ge Q-P
\end{equation}
Now, if $J=I_1$, $P(\tau)=\tau+g_3(\tau)$, $Q(\tau)
=P(\tau)+g_1(\tau)$ and there exists $\theta \in [0,1]$ such that
$T_0= \tau + g_3(\tau) + \theta g_1(\tau)$. Substituting in
\eqref{tre:venti} and using (\ref{F5}) we obtain
\begin{equation}\label{giuno}
g_1 (\tau) \le \frac{2\sqrt{2}}{c} \frac{f(\tau +g_3(\tau)+ \theta
g_1(\tau))}{f'(\tau +g_3(\tau)+\theta g_1(\tau))} \le
\frac{2D\sqrt{2}}{c} (\tau +g_3(\tau)+\theta g_1(\tau))^\mu
\end{equation}
In case $\mu \le 0$ we immediately obtain $g_1(\tau) = O(\tau)$,
hence we examine the case $\mu \in (0,1)$. Using the already known
equality $g_3(\tau)=O(\tau)$ and inequality $(x+y)^\mu \le x^\mu +
y^\mu$, there exist constants $K_1,K_2>0$ such that
\begin{equation}\label{stimagiuno}
\frac{g_1(\tau)}{\tau} \le  \frac{K_1}{\tau^{1-\mu}} + \frac{K_2
g_1(\tau)^\mu}{\tau}.
\end{equation}
Using a simple reasoning by contradiction, (\ref{stimagiuno})
implies $g_1(\tau)=O(\tau)$ as $\tau \to +\infty$.\par

If $J=I_1'$, $P(\tau)=\tau + (g_1+g_2+g_3)(\tau) +g_3'(\tau)$,
$Q(\tau)=P(\tau)+g'_1(\tau)$, $T_0=\tau + (g_1+g_2+g_3)(\tau)
+g_3'(\tau) + \theta g'_1(\tau)$, and substituting into
\eqref{tre:venti}
\begin{equation}
g_1'(\tau) \le \frac{2\sqrt{2}}{c} \frac{f(\tau+(g_1+g_2+g_3)
(\tau) +g_3'(\tau)+\theta g_1'(\tau))}{f'(\tau+(g_1+g_2+g_3)
(\tau) +g_3'(\tau)+\theta g_1'(\tau))} \label{eq:nuova}
\end{equation}
We will come back to this inequality later to prove
$g_1'(\tau)=O(\tau)$ as $\tau\to +\infty$. Indeed, by the same
argument as above, the only things that remain to show for this
purpose are $g_2(\tau)=O(\tau)$ and $g_3'(\tau)=O(\tau)$ as
$\tau\to +\infty$, and we are going to
prove these facts now.\\
\par
We consider an interval $J$ of type $2$ and again let $P(\tau) <
Q(\tau)$ denote its end points. Clearly $y(P)=1$ and $y(Q) = \mu
>1$ (or $y(Q^-)=+\infty$ in case that $z(Q)=0$. Indeed, what
follows works with any $\mu>0$). Again
$$
y' \ge 2 \sqrt{A(t)} y \qquad \text{and} \qquad \frac{y'}{y^2} \ge
\frac{1}{v} \qquad \text{a.e. on } J
$$
Fix $t \in[P,Q)$. Using $y(P)=1$, integration of the first
inequality on $[P,t]$ yields
\begin{equation}\label{tre:quindici}
y(t) \ge \exp\Big\{2 \int_P^t \sqrt{A(s)} ds\Big\} \qquad \forall
\ t \in [P,Q),
\end{equation}
while integrating the second one on $[t,Q-\varepsilon)$, for some
small $\varepsilon >0$, and proceeding as in (\ref{tre:a}) we have
\begin{equation}\label{tre:diciasette}
\frac{1}{y(t)} \ge \int_{t}^Q \frac{ds}{f(s)} \qquad \forall \ t
\in (P,Q).
\end{equation}
Thus, observing that
$$
2 \int_P^{t} {\sqrt{A(s)} ds} \ge \log \left( \frac{f(t)}{f(P)}
\right)^c
$$
we deduce from \eqref{tre:quindici}
$$
\frac{1}{y(t)} \le \left( \frac{f(P)}{f(t)} \right)^c.
$$
Finally, substituting into \eqref{tre:diciasette}
\begin{equation}\label{tre:diciannove}
1 \ge \left( \frac{f(t)}{f(P)} \right)^c\int_{t}^Q \frac{ds}{f(s)}
\ge (Q-t) \frac{1}{f(Q)} \left( \frac{f(t)}{f(P)} \right)^c \quad
\forall \ t \in (P,Q)
\end{equation}
Suppose now $J=I_2$ so that $g_2 (\tau)=Q(\tau)-P(\tau)$,
$$
\begin{array}{l}
P(\tau)=\tau+g_3(\tau)+g_1(\tau)\\
Q(\tau)= \tau+g_3(\tau)+g_1(\tau)+g_2(\tau)
\end{array}
$$
and since $t \in (P,Q)$, for some $\xi \in (0,1)$ we have
$$
\begin{array}{l}
t=\tau+(1-\xi)g_2(\tau) +g_1(\tau) +g_3(\tau)\\
Q-t=\xi g_2(\tau)
\end{array}
$$
Substituting into \eqref{tre:diciannove} yields,
\begin{equation}\label{gidue}
\sup_{\xi \in (0,1)} {\frac{\xi g_2(\tau)f(\tau+(1-\xi)
g_2(\tau)+g_1(\tau)+g_3(\tau))^c
}{f(\tau+g_2(\tau)+g_1(\tau)+g_3(\tau))f(\tau+g_1(\tau)+g_3(\tau))^c}}
\le 1
\end{equation}
Thus, setting $(g_1+g_3)(\tau)=k(\tau)$ since $g_1(\tau) =O(\tau)$
and $g_3(\tau)=O(\tau)$ as $\tau \to +\infty$, we have that
$k(\tau)=O(\tau)$ as $\tau \to +\infty$ and
$$
g_2 \in {\cal B}(f,k,c)
$$
and so $g_2(\tau) = O(\tau)$ as $\tau \to +\infty$. \\
\par
We can now deal with the case $J=I_3'$. We have already shown that
$g_1(\tau)+g_2(\tau)+g_3(\tau)=O(\tau)$ as $\tau \to \infty$. We
go back to \eqref{tre:quattordici} with
$J=I_3'=(P(\tau),Q(\tau))$: note that now
$$
P(\tau)=\tau +g_3(\tau)+g_1(\tau)+g_2(\tau) \quad , \quad
Q(\tau)=P(\tau)+ g_3' (\tau)
$$
where, obviously, $g_3'(\tau)=|I_3'|$. Since $t \in (P,Q)$, for
some $\xi \in (0,1)$ we have
$$
\begin{array}{l}
t = \tau + (1-\xi)g_3' (\tau) + (g_3+g_1+g_2)(\tau)\\
t-P=(1-\xi)g_3'(\tau)
\end{array}
$$
and substituting into \eqref{tre:quattordici}, since $t \in
(P,Q)$, is arbitrary we have
\begin{equation}\label{gitreprimo}
\sup_{\xi \in (0,1)} \frac{(1-\xi) g_3'(\tau) f(\tau+g_3' (\tau) +
(g_1+g_2+g_3)(\tau))^c}{f(\tau+(1-\xi)g_3'(\tau)+(g_1+g_2+g_3)(\tau))^{c+1}}
\le 1
\end{equation}
Thus, setting $h(\tau)= (g_1+g_2+g_3)(\tau)$, $h(\tau)= O (\tau)$
as $\tau \to +\infty$ and so we have $g_3' \in {\cal A}(f,h,c)$
therefore $g_3' (\tau) =O(\tau)$ as $\tau \to +\infty$. \\
\par
Coming back to inequality \eqref{eq:nuova}, we can now claim that
also $g_1'(\tau) = O(\tau)$ as $\tau \to +\infty$. \\
\par
The last case is $J=I_2'$ so that $g_2'(\tau)=Q(\tau)- P(\tau)$.
Now we have
$$
\begin{array}{l}
P(\tau)=\tau +(g_3+g_1+g_2+g_3'+g_1') (\tau)\\
Q(\tau)=P(\tau)+g_2'(\tau)
\end{array}
$$
and since $t\in(P,Q)$ there exists $\xi \in (0,1)$ such that
$$
\begin{array}{l}
t= \tau + (1-\xi) g_2'(\tau)+ (g_3+g_1+g_2+g_3'+g_1')(\tau)\\
Q(\tau)-t= \xi g_2'(\tau)
\end{array}
$$
Setting $k(\tau)=(g_3+g_1+g_2+g_3'+g_1')(\tau)$, we have already
proved that $k(\tau)=O(\tau)$ as $\tau \to +\infty$. Substituting
into \eqref{tre:diciannove} yields
\begin{equation}\label{gidueprimo}
\displaystyle \sup_{\xi \in (0,1)} \frac{\xi g_2'(\tau)
f(\tau+(1-\xi)g_2'(\tau)+k(\tau))^c}{f(\tau+g_2'(\tau)+k(\tau))f(\tau
+k(\tau))^c}\le 1
\end{equation}
Thus we have
$$
g_2' \in {\cal B}(f,k,c)
$$
therefore $g_2'(\tau)=O(\tau)$ as $\tau \to +\infty$, and this
shows that
$$
T_2(\tau) - T_1(\tau) \le T_2(\tau) - \tau =
(g_3+g_1+g_2+g_3'+g_1'+g_2')(\tau) = O(\tau)
$$
as $\tau \to +\infty$, so we have the first part of the theorem, that is, \eqref{tre:sei}.\\
\par
To conclude, we shall estimate the quantity
\begin{displaymath}
K = \limsup_{\tau \to +\infty} \frac{T_2(\tau) - \tau}{\tau}.
\end{displaymath}
Looking at the group of equations \eqref{gitre}, \eqref{giuno},
\eqref{gidue}, \eqref{gitreprimo}, \eqref{eq:nuova} and
\eqref{gidueprimo}, we first note that each of the functions
$g_i(\tau)$ and $g_i'(\tau)$ involved in the proof (shortly
$g(\tau)$) satisfies one of the following inequalities, for $\tau
\ge T$ and for some suitable function $h(\tau)$ which is known to
be $O(\tau)$:
\begin{align}
& \displaystyle \sup_{\xi \in (0,1)} \frac{(1-\xi) g(\tau)
f(\tau+g(\tau)+h(\tau))^c}{f(\tau+(1-\xi)g(\tau)+h(\tau))^{c+1}}
\le 1 &  \quad {\rm for } \ g_3 \ {\rm and } \ g_3' & &
\label{casodeltre}
\\[0.2cm]
& \displaystyle g(\tau) \le \frac{2\sqrt{2}}{c}\frac{f(\tau+
h(\tau)+\theta g(\tau))}{f'(\tau + h(\tau)+\theta g(\tau))} &
\quad {\rm for } \ g_1 \ {\rm and } \ g_1' & & \label{casodelluno}
\\[0.2cm]
& \displaystyle \sup_{\xi \in (0,1)} \frac{\xi g(\tau)
f(\tau+(1-\xi)g(\tau)+h(\tau))^c}{f(\tau+g(\tau)+h(\tau))\cdot
f(\tau+h(\tau))^c} \le 1 & \quad {\rm for } \ g_2 \ {\rm and } \
g_2' & & \label{casodeldue}
\end{align}\\
For the sake of simplicity, we perform computations in case
$$
f(t) = \Lambda \exp\{at^\alpha\}, \qquad a, \Lambda, \alpha > 0
\quad
$$
(note that $f$ satisfy property $(P)$ for every $c>1$). We shall
determine $K$ by computing, in each of the tree cases above,
$$
K_j = \limsup_{\tau\to +\infty} \frac{g(\tau)}{\tau} \qquad
$$
(the index $j$ corresponds to the cases satisfied by $g_j$ and
$g_j'$), and then summing the terms "inductively" following the
changes of the known function $h$ case by case. For this purpose
let
$$ H \ge \limsup_{\tau\to +\infty}\frac{h(\tau)}{\tau}
$$\par
Consider at first inequality \eqref{casodelluno}: we immediately
find that, for this choice of $f$,
$$
\frac{g(\tau)}{\tau} \le
\frac{2\sqrt{2}}{c}\frac{1}{\tau}\frac{1}{a\alpha
(\tau+h(\tau)+\theta g(\tau))^{\alpha-1}} \le
\frac{2\sqrt{2}}{ca\alpha} \frac{\big(1 + \frac{h(\tau)}{\tau} +
\frac{g(\tau)}{\tau}\big)^{1-\alpha}}{\tau^\alpha}
$$
We claim that $K_1=0$. Indeed, suppose by contradiction that there
exists a divergent sequence $\{\tau_n\}$ such that
$g(\tau_n)/\tau_n \rightarrow K_1>0$. Then, evaluating the above
inequality along $\{\tau_n\}$ and passing to the limit we reach
$$
0 < K_1 \le 0 \qquad \mathrm{a \ contradiction.}
$$
We now focus our attention on \eqref{casodeltre}. By an algebraic
manipulation
$$
g(\tau) \le \frac{1}{1-\xi }
\frac{f(\tau+(1-\xi)g(\tau)+h(\tau))^{c+1}}{f(\tau+ g(\tau) +
h(\tau))^c} \quad \forall \ \xi \in (0,1) \ .
$$
Due to the form of $f$, better estimates can be obtained choosing
$\xi$ near $1$. For $\tau >1$, we choose $\xi=(\tau-1)/\tau$. For
the ease of notation let $x(\tau)=g(\tau)/\tau$, so that $x(\tau)$
is bounded on $[T,+\infty)$ because $f$ satisfies property $(P)$.
With this choice of $\xi$ we have
\begin{equation}\label{casoterzo}
x(\tau) \le \frac{f(\tau+x(\tau)+h(\tau))^{c+1}}{f(\tau+\tau
x(\tau) + h(\tau))^c},
\end{equation}
thus substituting
$$
x(\tau) \le \Lambda \exp\Big\{a\tau^\alpha \Big[
(c+1)\Big(1+\frac{x(\tau)}{\tau} + \frac{h(\tau)}{\tau}
\Big)^\alpha -c\Big(1+x(\tau) + \frac{h(\tau)}{\tau} \Big)^\alpha
\Big] \Big\}
$$
Suppose now that $K_3>0$, and evaluate this inequality along a
sequence $\{\tau_n\}$ such that $x(\tau_n) \to K_3$. Choose $0
<\delta <K_3 $, and let $n$ be large enough that the following
inequalities hold:
$$
x(\tau_n) > K_3 -\delta \quad, \quad \frac{x(\tau_n)}{\tau_n}
<\delta
$$
This yields:
\begin{equation}\label{stimatre}
x(\tau_n) \le \Lambda \exp\Big\{a\tau_n^\alpha \Big[ (c+1)\Big(1+
\delta + \frac{h(\tau_n)}{\tau_n} \Big)^\alpha -c\Big(1+ K_3
-\delta + \frac{h(\tau_n)}{\tau_n} \Big)^\alpha \Big] \Big\}
\end{equation}
Suppose now that $K_3$ satisfies
\begin{equation}\label{stimakappatre}
\max_{\mu \in [0,H]} \Big\{(c+1) (1+\mu)^\alpha -
c(1+K_3+\mu)^\alpha \Big\} <0,
\end{equation}
and compare it with \eqref{stimatre}. We can say that, by
continuity, there exists a small $\delta >0$ such that the
expression between square brackets is strictly less than $0$.
Letting now $\tau_n$ go to infinity in \eqref{stimatre} we deduce
$0<K_3 \le 0$, a contradiction. Note that \eqref{stimakappatre}
holds if and only if
$$
(c+1) - c\Big(\frac{K_3}{\mu+1} +1\Big)^\alpha <0 \qquad \forall \
\mu \in [0,H],
$$
that is,
$$ K_3 > \Big[\Big(\frac{c+1}{c}\Big)^{\frac{1}{\alpha}}
-1\Big](1+H).
$$
Hence, if $K_3>0$, we necessarily have
\begin{equation}\label{kappatre} K_3 \le
\Big[\Big(\frac{c+1}{c}\Big)^{\frac{1}{\alpha}} -1\Big](1+H).
\end{equation}
\par
The same technique can be exploited when dealing with
\eqref{casodeldue}: from
\begin{equation}\label{casosecondo}
g(\tau) \le \frac{1}{\xi} \frac{f(\tau+g(\tau)+h(\tau))\cdot
f(\tau + h(\tau))^c}{f(\tau+ (1-\xi)g(\tau)+h(\tau))^c} \qquad
\forall \ \xi \in (0,1)
\end{equation}
we deduce that it is better to choose $\xi$ near 0, so we set $\xi
=1/\tau$ and we obtain, with the same notations,
$$
x(\tau) \le \frac{f(\tau + \tau x(\tau) + h(\tau))\cdot f(\tau +
h(\tau))^c}{f(\tau+(\tau-1)x(\tau)+h(\tau))^c}
$$
Thus
$$
\begin{array}{lcl}
x(\tau) & \le & \Lambda \exp\Big\{a\tau^\alpha \Big[
\Big(1+x(\tau)+\frac{h(\tau)}{\tau}\Big)^\alpha + \\[0.4cm]
 & & + c\Big( 1+ \frac{h(\tau)}{\tau} \Big)^\alpha -c\Big( 1 +
\frac{\tau-1}{\tau} x(\tau) + \frac{h(\tau)}{\tau}\Big)^\alpha
\Big] \Big\}
\end{array}
$$
Next, if $K_2>0$ we choose a sequence $\{\tau_n\}$ realizing $K_2$
and we consider $n$ sufficiently large that
$$
\frac{(\tau_n-1)}{\tau_n} > (1-\delta) \quad , \quad K_2-\delta <
x(\tau_n)<K_2+\delta
$$
obtaining the estimate
\begin{equation}
\begin{array}{lcl}
x(\tau_n) & \le & \Lambda \exp\Big\{a\tau_n^\alpha \cdot \Big[
\Big(1+ (K_2+\delta)+\frac{h(\tau_n)}{\tau_n}\Big)^\alpha +
\\[0.4cm]
& & + c\Big( 1+ \frac{h(\tau_n)}{\tau_n} \Big)^\alpha -c\Big( 1 +
(1-\delta)(K_2-\delta) + \frac{h(\tau_n)}{\tau_n}\Big)^\alpha
\Big] \Big\}\label{stimadue}
\end{array}
\end{equation}
Now, if $K_2$ satisfies
\begin{equation}\label{stimakappadue}
\max_{\mu \in [0,H]} \Big\{(1+K_2 +\mu)^\alpha +c(1+\mu)^\alpha
-c(1+K_2+\mu)^\alpha \Big\}<0,
\end{equation}
we reach a contradiction proceeding as in the previous case.
Similarly to what we did above this yields the bound
\begin{equation}\label{kappadue}
K_2 \le \Big[\Big(\frac{c}{c-1}\Big)^{\frac{1}{\alpha}}
-1\Big](1+H)
\end{equation}\par
To simplify the writing we now set
$$
W = \Big[\Big(\frac{c+1}{c}\Big)^{\frac{1}{\alpha}} -1\Big] \quad
, \quad Z = \Big[\Big(\frac{c}{c-1}\Big)^{\frac{1}{\alpha}}
-1\Big]
$$
To estimate $g_3(\tau)/\tau$, we shall use \eqref{kappatre} and,
from \eqref{gitre}, we deduce $h(\tau) \equiv 0$ and thus $H=0$.
Therefore, we get
$$
K_3 \le W.
$$
We have already shown that $K_1=0$. Next, to estimate
$g_2(\tau)/\tau$ we shall consider \eqref{kappadue}. By
\eqref{gidue} $h(\tau) = g_3(\tau)+g_1(\tau)$, so we can use for
$H$ the sum $W+0 = W$, hence
$$
K_2 \le Z(1+W).
$$
Proceeding along the same lines we obtain the estimates
$$
\begin{array}{l}
K_3' \le W\big(1+ W + Z(1+W)\big); \\[0.2cm]
K_1' = 0;\\[0.2cm]
K_2' \le Z\Big(1+ W + Z(1+W) + W\big(1+W+Z(1+W)\big)\Big).
\end{array}
$$
Summing up the $K_j$ and the $K_j'$, we obtain the surprisingly
simple expression
$$
K\le \sum_{j=1}^3 (K_j + K_j') = \big(W+1 \big)^2\big(Z+1 \big)^2
- 1 = \Big(\frac{c+1}{c-1}\Big)^{\frac{2}{\alpha}} -1 \quad .
$$
Thus we eventually have
\begin{equation}\label{stimalimsup}
\limsup_{\tau\to +\infty} \frac{T_2(\tau)}{\tau} \le
\Big(\frac{c+1}{c-1}\Big)^{\frac{2}{\alpha}}
\end{equation}
With few modifications in the computations, it can be seen that,
considering $f(t)=\Lambda \exp[at^\alpha \log^\beta t]$ instead of
the above, the value of the constant $K$ does not change.
$\blacktriangle$
\begin{Remark}\label{ft}
\emph{Since $f(t)= \Lambda \exp\{at^\alpha \log^\beta t\}$
satisfies property $(P)$ for every $c>1$, in this case conditions
(\ref{A5}) and (\ref{A4}) may be replaced by
\begin{equation}\tag{A3 + A4}\label{A3+A4}
\sqrt{A(t)} \ge c\Big(\frac{a \alpha}{2}\Big) t^{\alpha-1}
\log^\beta t \qquad \text{a.e. on } [T,+\infty), \text{ for some }
c>1.
\end{equation}
Note that the right hand side of the above expression is
asymptotic to $c\sqrt{\widetilde{\chi}_f(r)}$, hence to
$c\sqrt{\chi_f(t)}$ by (\ref{asinttt}). Indeed, choose
$\varepsilon>0$ such that $\overline{c}=c-2\varepsilon >1$. Then,
there exists $T>0$ such that, on $[T,+\infty)$,
$$
c \Big(\frac{a \alpha}{2}\Big) t^{\alpha-1} \log^\beta t \ge
(c-\varepsilon) \sqrt{\chi_f(t)} \ge
(c-2\varepsilon)\sqrt{\widetilde{\chi}_f(t)}=
\overline{c}\sqrt{\widetilde{\chi}_f(t)},
$$
hence (\ref{A4}) is satisfied with $\overline{c}$. Moreover,
$(iv)$ of Proposition \ref{condsuff} and $\chi_f(t) \ge \chi(t)$
ensure condition (\ref{A5}). Applying Theorem \ref{teo:princ} with
$\overline{c}$ yields
$$
\limsup_{\tau \rightarrow +\infty} \frac{T_2(\tau)}{\tau} =
\Big(\frac{\overline{c}+1}{\overline{c}-1}\Big)^{\frac{2}{\alpha}}
$$
Since $\varepsilon$ is arbitrary, we conclude the validity of
(\ref{limsupa}) under assumption (\ref{A3+A4}).}
\end{Remark}

\begin{Remark}
\emph{One might ask if varying the choice of the level sets in
\eqref{livelli} one could obtain better estimates. It is not hard
to see that, for every choice of the level, \eqref{stimalimsup}
does not change.}
\end{Remark}
%
%
%
%
%
%
\section{Geometric applications}\label{geometrici}
This section is devoted to the proofs of the geometric
applications given in the Introduction, which follow from the
results of sections \ref{existencefirstzero} and
\ref{oscillestimates}. The core are Theorems \ref{spettro} and
\ref{principale}, where the Cauchy problems (\ref{cauchy2}),
(\ref{cauchy2inf}) appear in order to obtain suitable radial test
functions which yield estimates for the Rayleigh quotients of $L$
and $\Delta$ respectively. An almost direct use of Theorem
\ref{spettro} proves Theorems \ref{main}, \ref{Gaussmap} and
\ref{Yamabe}, while Theorem \ref{superfici} requires some special
attention and further work.

\subsection{The index of $\Delta + a(x)$: proof of Theorem \ref{spettro}}

Choose $v(t)=\mathrm{Vol}(\partial B_t)$. From Proposition
\ref{proprietavol} it follows that the spherical mean $A(t)$
belongs to $L^\infty_\mathrm{loc}([0,+\infty))$, the validity of
\eqref{V} and the existence of a locally Lipschitz solution of
\eqref{cauchy2} whose zeros (if any) are isolated (Theorems
\ref{esist} and \ref{isolati} of the Appendix). Consider problem
(\ref{cauchy2}), and note that, by the coarea formula,
$$
0 < \int_0^{R_0} A(s)v(s) ds = \int_0^{R_0} \Big( \int_{\partial
B_s} a\Big) ds = \int_{B_{R_0}} a.
$$
By Corollary \ref{cor:zero}, assumption $(i)$ guarantees the
existence of a first zero of every locally Lipschitz solution
$z(t)$, whereas Corollary \ref{corforte} implies that assumption
$(ii)$ forces $z(t)$ to be oscillatory. Note that a different
choice of $R$ in assumption $(ii)$ does not affect the value of
the ``limsup".\par
We now consider case $(i)$: choose a locally Lipschitz solution
$z(t)$ of (\ref{cauchy2}), and denote with $T$ its first zero.
Define
$$
\psi(x) = z(r(x))
$$
so that
$$
\psi \in \mathrm{Lip}(\overline{B}_T) \quad , \quad \psi \equiv 0
\quad \text{on }
\partial B_T \quad , \quad \nabla \psi (x) = z'(r(x))\nabla r(x)
\quad \text{a.e. on } M
$$
and fix $0<\varepsilon<T$. Then, using the coarea formula, Gauss
lemma and (\ref{cauchy2}) we obtain
$$
\begin{array}{l}
\displaystyle  \int_{B_T\backslash B_\varepsilon} |\nabla \psi|^2
- a(x)\psi^2 =  \int_{B_T\backslash B_\varepsilon}
|\nabla \psi|^2 - A(r)\psi^2 \\[0.4cm]
 = \displaystyle \int_\varepsilon^T (z'(r))^2v(r)dr -
\int_\varepsilon^T A(r)z^2(r)v(r) dr \\[0.4cm]
 =  \displaystyle
-z(\varepsilon)z'(\varepsilon)v(\varepsilon) - \int_\varepsilon^T
z(r)\big[(v(r)z'(r))'+A(r)v(r)z(r)\big] =
-z(\varepsilon)z'(\varepsilon)v(\varepsilon)
\end{array}
$$
and letting $\varepsilon\downarrow 0^+$ we deduce
$$
\int_{B_T} |\nabla \psi|^2 -a(x)\psi^2 \le 0.
$$
By Rayleigh characterization of eigenvalues and by domain
monotonicity we conclude $\lambda_1^L(M)<0$.\par
Suppose now we are in case $(ii)$, and assume by contradiction
that there exists $R>0$ such that
\begin{equation}\label{fish}
\lambda_1^L(M\backslash B_R) \ge 0.
\end{equation}
As already stressed in the Introduction, by a result of
Fisher-Colbrie \cite{FC}, if the index of $L$ is finite then
(\ref{fish}) holds for a sufficiently large $R$.\par
In our assumptions, every locally Lipschitz solution $z(t)$ of
(\ref{cauchy2inf}) is oscillatory. Let $T_1<T_2$ be two
consecutive zeros of $z(t)$ strictly after $R$. Define $\psi(x)=
z(r(x))$ in the annular region $B_{T_2}\backslash B_{T_1}$, and
$\psi(x) \equiv 0$ in the rest of $M$. Then $\psi \in
\mathrm{Lip}_0(M)$ with support contained in $M\backslash B_R$.
Proceeding as in the previous case, we obtain
$$
\int_{B_{T_2}\backslash B_{T_1}} |\nabla \psi|^2 - a(x)\psi^2 \le
0,
$$
hence, by strict domain monotonicity, $\lambda_1(M\backslash B_R)
<0$, contradicting (\ref{fish}).\par
Let us finally consider case $(iii)$. By Remark \ref{ft} and
Theorem \ref{teo:princ}, (\ref{cauchy2}) is oscillatory, thus $L$
is unstable at infinity. In particular, the index of $L$ is
infinite. Note that (\ref{indexgrowth}) is equivalent to prove
that
$$
\liminf_{r\rightarrow +\infty} \frac{\mathrm{ind}_L(B_r)}{\log r}
\ge \frac{1}{\log K} \quad , \quad \mathrm{with } \quad K=
\Big(\frac{c+1}{c-1}\Big)^{\frac{2}{\alpha}}.
$$
Fix $\varepsilon>0$. Then, by Theorem \ref{teo:princ} there exists
$T=T(\varepsilon)$ such that on $[T,+\infty)$
$$
\frac{T_2(r)}{r} \le K_\varepsilon =
\Big(\frac{c+1}{c-1}\Big)^{\frac{2}{\alpha}} +\varepsilon.
$$
Proceeding as above, on $M\backslash B_r$ we can find a radial
function $\psi_1(x)$, with support contained in $B_{K_\varepsilon
r}$, which makes the Rayleigh quotient non positive. Starting from
$T_2(r)$, the second zero after $T_2(r)$ is attained before
$K_\varepsilon T_2(r) \le K_\varepsilon^2 r$, and we can construct
a new Lipschitz radial function $\psi_2(x)$ which makes the
Rayleigh quotient non positive. Moreover, the support of $\psi_2$
is disjoint from that of $\psi_1$. In conclusion, the index of $L$
grows at least by $1$ when the radius is multiplied by
$K_\varepsilon$, hence
$$
\mathrm{ind}_L(B_r) \ge \mathrm{ind}_L(B_T) +
\Big\lfloor\log_{K_\varepsilon} \Big(\frac{r}{T}\Big)\Big\rfloor,
$$
where $\lfloor s \rfloor$ denotes the floor of $s$. Therefore we
have
\begin{equation}\label{Keps}
\liminf_{r\rightarrow +\infty}
\frac{\mathrm{ind}_L(B_r)}{\log_{K_\varepsilon} r} \ge 1 \qquad
\forall \ \varepsilon >0.
\end{equation}
From the change of base theorem, for every $u,v > 1$, $r>0$
\begin{equation}\label{delop}
\frac{\log_u r}{\log_v r} = \log_u v = \frac{\log v}{\log u},
\end{equation}
so that
\begin{equation}
\liminf_{r\rightarrow +\infty} \frac{\mathrm{ind}_L(B_r)}{\log r}
 \ge \frac{1}{\log K_\varepsilon} \qquad \forall \ \varepsilon
>0.
\end{equation}
Letting $\varepsilon \rightarrow 0$ yields the desired conclusion.
$\blacktriangle$\\
\subsection{Tangent envelopes: proof of Theorem \ref{main}}

We briefly recall some well known facts. Suppose we are given an
isometrically immersed hypersurface
$$
\varphi : M^m \longrightarrow N^{m+1},
$$
where $N$ is orientable. We fix the index notation $i,j,k,t \in \{1,\ldots,m\}$, and we choose a local Darboux
frame $\{e_i,\nu\}$. Let $R, \mathrm{Ricc},s$ (resp $\overline{R}, \overline{\mathrm{Ricc}},\overline{s}$) be
the curvature tensor, the Ricci tensor and the scalar curvature of $M$ (resp. $N$), denote with $II= (h_{ij})$
the second fundamental form of the immersion, with $|II|^2$ the square of its Hilbert-Schmidt norm and with $H =
m^{-1}h_{ii} \nu$ the mean curvature vector. Tracing twice the Gauss equations
\begin{equation}\label{gaussequations}
R_{ijkt}= \overline{R}_{ijkt} + h_{ik}h_{jt} -h_{it} h_{jk}
\end{equation}
we get
\begin{equation}\label{scalareeq}
s= \overline{s} -2 \overline{\mathrm{Ricc}}(\nu,\nu) + m^2|H|^2 -
|II|^2.
\end{equation}
Moreover, we recall the Codazzi-Mainardi equation
\begin{equation}\label{Codazzi}
h_{ijk} - h_{ikj} = \overline{R}^{m+1}_{ijk},
\end{equation}
where $(h_{ijk})$ are the components of the covariant derivative $\nabla II$. A minimal immersion $\varphi$ is
characterized by $H\equiv 0$, which implies that $\varphi$ is a stationary point for the volume functional on
every relatively compact domain with smooth boundary in $M$. It is known that if, for example,
$N=\mathbb{R}^{m+1}$, a minimal hypersurface cannot be
compact and, by (\ref{scalareeq}), $s(x) = -|II|^2 \le 0$.\\
We say that $\varphi$ is stable if it locally minimizes the volume
functional up to second order, and unstable otherwise.
Analytically the condition of stability is expressed by
$$
\int_M |\nabla \psi|^2 -
\big(|II|^2+\overline{\mathrm{Ricc}}(\nu,\nu)\big)\psi^2 \ge 0
\qquad \forall\ \ \psi \in C^\infty_0(M)
$$
and it is equivalent to the fact that the Schr\"odinger operator
$L= \Delta + |II|^2 + \overline{\mathrm{Ricc}}(\nu,\nu)$ satisfies
$\lambda_1^L(M)\ge 0$. Observe that, if $N$ is Ricci flat (for
example, $N=\mathbb{R}^{m+1}$), using \eqref{scalareeq} we get
$$
L = \Delta + |II|^2 = \Delta - s(x).
$$
The strategy of the proof of Theorem \ref{main} is to proceed by
contradiction. First we prove that, if (\ref{tuttoRm}) fails, $M$
is stable at infinity, i.e. $\lambda_1^L(M\backslash \Omega) \ge
0$ for the chosen compact domain $\Omega$; then, we contradict
this fact using Theorem \ref{spettro} under assumptions
(\ref{volumenoninl1}) or (\ref{volume}), (\ref{mediasferica}),
(\ref{parametri}).\\
\par
\noindent \textbf{Proof (of Theorem \ref{main}). } We reason by
contradiction and, without loss of generality, we can assume that
the origin $o$ of $\mathbb{R}^{m+1}$ belongs to
$$
\mathbb{R}^{m+1}\backslash \bigcup_{x\in M\backslash \Omega} T_xM.
$$
Consider on $M\backslash \Omega$ a local normal unit vector field
$\nu$ and define the local vector field $X= \langle
\varphi,\nu\rangle \nu$, where $\langle,\rangle$ denotes the
canonical metric on $\mathbb{R}^{m+1}$. For every point $x$ in the
domain of $X$ we have $X_x\not \equiv 0$ since otherwise
$\varphi(x)$ would be orthogonal to $\nu(x)$ and thus $T_xM$ would
contain the origin $o$. Moreover, under a change of Darboux frame
the value of $X$ does not change, hence it provides a globally
defined, nowhere vanishing normal vector field, proving that
$M\backslash \Omega$ is orientable. Define $u(x)= \langle
\varphi(x),\nu(x)\rangle \neq 0$, $u \in C^\infty(M\backslash
\Omega)$. Possibly inverting the orientation on connected
components, we can suppose $u>0$ on $M\backslash \Omega$. A simple
computation using minimality of $\varphi$ and Codazzi equation
\eqref{Codazzi} for $N=\mathbb{R}^{m+1}$ shows that $u$ is a
positive solution of
$$
\Delta u - s(x) u = 0 \qquad \mathrm{on \ } M\backslash \Omega,
$$
By the result of Fisher-Colbrie and Schoen \cite{FCS} it follows
that $L = \Delta - s(x)$ has non negative spectral radius
$\lambda_1^L(M\backslash \Omega)$, hence $M$ is stable at infinity.\\
To contradict this latter result, we choose $A(r)= - S(r)$ and we
use Theorem \ref{spettro}, case $(ii)$. A contradiction is
immediate in case of \eqref{volumenoninl1}, while if we assume
$(\mathrm{Vol}(\partial B_r))^{-1} \in L^1(+\infty)$ we can apply
Proposition \ref{condsuff}, item $(iv)$: indeed, under assumptions
(\ref{volume}), \eqref{mediasferica} and (\ref{parametri}),
observing that
$$
\frac{d}{ds}(-s^{1-\alpha}\exp\{-s^\alpha\}) \le \widetilde{C}
\exp\{-s^\alpha\} \qquad \text{for } s\ge 1,
$$
for some $\widetilde{C}>0$, there exist positive constants $D$ and
$H$ such that
$$
\begin{array}{lcl}
\displaystyle \liminf_{r\rightarrow +\infty} \frac{\int_R^r
\sqrt{A(s)}ds}{-\frac{1}{2} \log \int_r^{+\infty}
\frac{ds}{\mathrm{Vol}(\partial B_s)}} & \ge & \displaystyle
\liminf_{r\rightarrow +\infty} \frac{\int_R^r \sqrt{C} s^{-\mu/2}
ds}{-\frac{1}{2} \log \int_r^{+\infty}
\frac{\exp\{-s^\alpha\}}{\Lambda }ds } \\[0.6cm]
& \ge & \displaystyle  \liminf_{r\rightarrow +\infty} \Big(D
r^{1-\frac{\mu}{2}-\alpha} \log^{-H} r\Big) = +\infty.
\end{array}
$$
Proposition \ref{condsuff} item $(iv)$ implies (\ref{infinite}),
so that Theorem \ref{spettro} case $(ii)$ contradicts the
stability at infinity of $L$.
$\blacktriangle$\\
\begin{Remark}
\emph{Obviously, when $\Omega=\emptyset$ there is a version of the
above theorem in finite form, which is based on case $(i)$ of
Theorem \ref{spettro}. We have preferred not to make the
proposition too cumbersome, in order to better appreciate the
result itself. Nevertheless, even this case seems interesting:
inequality (\ref{finite}) implies that a strongly negative scalar
curvature on a compact set spreads the tangent hyperplanes
everywhere on $\mathbb{R}^{m+1}$, independently of the behavior of
the curvature outside the compact.}\\
\end{Remark}
\subsection{The Gauss map: proof of Theorem \ref{Gaussmap}}
The proof follows the same lines of Theorem \ref{main}, and we
maintain the same notations. We fix an equator $E$ and we reason
by contradiction: assume that there exist a sufficiently large
geodesic ball $B_R$ such that, outside $B_R$, $\nu$ does not meet
$E$. In other words $\nu (M \backslash B_R)$ is contained in the
open spherical cups determined by $E$. Indicating with $w\in
\mathbb{S}^m$ one of the two focal points of $E$, we can say that
$\langle w,\nu(x)\rangle \neq 0$ for every $x\in M\backslash B_R$,
where $\langle,\rangle$ stands for the scalar product of unit
vectors in $\mathbb{S}^m\subset \mathbb{R}^{m+1}$. Then, the
normal vector field $X= \langle w,\nu\rangle \nu$ is globally
defined and nowhere vanishing on $M\backslash B_R$, proving that
$M\backslash B_R$ is orientable. Therefore, the Gauss map is
globally defined on $M\backslash B_R$. Let $\mathcal{C}$ be one of
the (finitely many) connected components of $M\backslash B_R$;
then, $\nu(\mathcal{C})$ is contained in only one of the open
spherical caps determined by $E$. Up to replacing $w$ with $-w$,
we can suppose $u=\langle w,\nu \rangle >0$ on $\mathcal{C}$.
Proceeding in the same way for every connected component, we can
construct a positive function $u$ on $M\backslash B_R$. By a
standard calculation $u$ satisfies
\begin{equation} \label{due:quattordici}
\Delta u = -h_{ikk}\langle e_i, w\rangle  - |II|^2u \qquad {\rm
on} \ M \backslash B_R.
\end{equation}
Using Schwarz symmetry, Codazzi equation \eqref{Codazzi} and minimality we deduce
$$
h_{ikk}=h_{kik}= h_{kki} = 0,
$$
hence $\Delta u + |II|^2 u = 0$. From \eqref{scalareeq} we get
\begin{equation}\label{equazstabile2}
\Delta u -s(x)u = 0
\end{equation}
In particular, \eqref{equazstabile2} implies
$\lambda_1^L(M\backslash B_R)\ge 0$, where $L=\Delta -s(x)$.
Observe that $s(x)= -|II|^2 \le 0$, so that its spherical mean
$S(r)$ is non positive. As in the proof of Theorem \ref{main}, the
assumptions imply case $(ii)$ of Theorem \ref{spettro}, and this
contradicts $\lambda_1^L(M\backslash B_R)\ge 0$.
$\blacktriangle$\\
\subsection{The Yamabe Problem: proof of Theorem \ref{Yamabe}}
Applying Theorem \ref{spettro} to the operator $L= \Delta -
\frac{1}{c_m}s(x)$ we obtain $\lambda_1^L(M)<0$. Hence, the
conclusion follows from Theorems $2.4$ and $2.1$ of \cite{PRS},
with the observation after Theorem $2.3$ therein. $\blacktriangle$
\begin{Remark}
\emph{We can state an alternative version at infinity of condition
(\ref{yam}) via Proposition \ref{condsuff}, $(iv)$. This reads as
follows. Suppose that
$$
S(r) \le -\frac{H}{r^\beta} \qquad \mathrm{for} \ r
\gg 1 \ \mathit{and \ some \ } H>0 , \ \beta \le 2.
$$
Then, condition
\begin{equation}
\sqrt{H\frac{m-2}{m-1}} > \left\{\begin{array}{ll} \displaystyle
\liminf_{r\rightarrow +\infty} \Big( \frac{\frac{\beta}{2}
-1}{r^{-\beta/2+1}} \log \int_r^{+\infty}
\frac{ds}{\mathrm{Vol}(\partial B_s)}\Big) & \quad
\mathit{if \ } \beta < 2 \\[0.6cm]
\displaystyle \liminf_{r\rightarrow +\infty} \Big( -\frac{1}{\log
r} \log \int_r^{+\infty} \frac{ds}{\mathrm{Vol}(\partial B_s)}
\Big) & \quad \mathit{if \ } \beta = 2
\end{array}\right.
\end{equation}
implies the existence of the desired conformal deformation.}\\
\end{Remark}
\subsection{Minimal surfaces: proof of Theorems \ref{superfici} and \ref{corpiatto}}
We will obtain both the results as easy consequences of the next two lemmas, the first of which is a somewhat
modified version of a result of Colding and Minicozzi \cite{CM}. We adopt the notations of Theorems \ref{main}
and \ref{Gaussmap}.
\begin{Lemma}\label{parab}
Let $\varphi:M^2\rightarrow N^3$ be a simply connected, minimally
immersed surface in an ambient $3$-manifold. Assume that the Ricci
tensor of $N$ satisfies
\begin{equation}\label{riccimagzero}
\overline{\mathrm{Ricc}} \ge 0.
\end{equation}
Suppose that $M$ has a pole $o$, and let $L$ be the stability
operator. If $\lambda_1^L(M\backslash \Omega) \ge 0$ for some
compact set $\Omega$, then there exists a constant $C>0$ such that
$$
\mathrm{Vol}(B_R) \le CR^2 \qquad \forall \ R \ge 0.
$$
\end{Lemma}
\noindent\textbf{Proof. } Let $K$ be the sectional curvature of
$M$. Since for surfaces $s(x)=2K$, using \eqref{riccimagzero} in
\eqref{scalareeq} yields
$$
\begin{array}{lcl}
2K & = & \displaystyle \overline{\mathrm{Ricc}}(e_1,e_1) + \overline{\mathrm{Ricc}}(e_2,e_2)
- \overline{\mathrm{Ricc}}(\nu,\nu) - |II|^2 \\[0.2cm]
& \ge & - \overline{\mathrm{Ricc}}(\nu,\nu) -|II|^2
\end{array}
$$
hence the Rayleigh quotient for the stability operator do not
exceed that for $\overline{L}=\Delta - 2K$. It follows that, for
every subset $D\subset M$, we have inequality
\begin{equation}\label{stimaL}
\lambda_1^L(D) \le
\lambda_1^{\overline{L}}(D),
\end{equation}
thus by the assumptions $\lambda_1^{\overline{L}}(M\backslash
B_{R_0})\ge 0$ for some $R_0$ sufficiently large that $\Omega\subset B_{R_0}$.\\
Since $M$ is simply connected and has a pole, the geodesic spheres
centered at $o$ are smooth and the geodesic balls are
diffeomorphic to Euclidean ones. By Gauss-Bonnet theorem together
with the first variation formula, we have
\begin{equation}\label{gb}
\int_{B_r} K = 2\pi - l'(r),
\end{equation}
where $l(r)$ is the length of $\partial B_r$ (another way to
derive this formula can be found in \cite{PRS1}, page $238$).
Denote with $K(r) = \int_{B_r} K$, and observe that, by the coarea
formula $K'(r) =
\int_{\partial B_r} K$.\\
By the stability of $\overline{L}$, for every $\psi \in
\mathrm{Lip}_0(M\backslash B_{R_0})$ we have
\begin{equation}\label{stab}
\int_{M\backslash B_{R_0}} |\nabla \psi|^2 + 2 \int_{M\backslash
B_{R_0}} K \psi^2 \ge 0.
\end{equation}
Fix $R>R_0+2$ and choose $\psi(x) = f(r(x))$, where
$$
f(t) = \left\{\begin{array}{ll}
0 & \mathrm{if \ } t \le R_0 \\[0.1cm]
\displaystyle t-R_0 & \mathrm{if \ } t \in [R_0,R_0+1] \\[0.1cm]
\displaystyle \frac{R-t}{R-R_0-1} & \mathrm{if \ } t \in [R_0+1,R]
\\[0.1cm]
0 & \mathrm{if \ } t \ge R
\end{array}\right.
$$
Then, using (\ref{gb}) into (\ref{stab}) and integrating by parts,
by the properties of $f$ we have
$$
0 \le \int_{R_0}^R (f'(r))^2l(r) dr  + 2 \int_{R_0}^R
l'(r)(f^2(r))'dr
$$
Inserting the explicit expression of $f$ we obtain
$$
\begin{array}{l}
\displaystyle 0 \le \mathrm{Vol}(B_{R_0+1}) -
\mathrm{Vol}(B_{R_0}) + \frac{\mathrm{Vol}(B_R) -
\mathrm{Vol}(B_{R_0+1})}{(R-R_0-1)^2} +
4 l(R_0+1) + \\[0.3cm]
\displaystyle -
4\big(\mathrm{Vol}(B_{R_0+1})-\mathrm{Vol}(B_{R_0})\big) +
\frac{4l(R_0+1)}{R-R_0-1} -
\frac{4\big(\mathrm{Vol}(B_R)-\mathrm{Vol}(B_{R_0+1})\big)}{(R-R_0-1)^2}
\end{array}
$$
Therefore, there exists a constant $C=C(R_0)$ depending on the
geometry of $B_{R_0+1}$ such that, for every $R > R_0+2$,
$$
\frac{3\big(\mathrm{Vol}(B_R)-\mathrm{Vol}(B_{R_0+1})\big)}{(R-R_0-1)^2}
\le C(R_0)
$$
hence,
$$
\mathrm{Vol}(B_R) \le \mathrm{Vol}(B_{R_0+1}) + \frac{C(R_0)}{3}
(R-R_0-1)^2 \le \widetilde{C(R_0)}R^2.
$$
Since near $o$ the geometry of $M$ is ``nearly" Euclidean, up to
enlarging the constant the same estimate holds on all of $M$, and
this
concludes the proof. $\blacktriangle$\\
\par
\begin{Remark}\label{lift}
\emph{Note that, in case $\Omega=\emptyset$ and $N=\mathbb{R}^3$,
we recover Colding and Minicozzi theorem, for which the
simply-connectedness assumption is unnecessary: in fact, we can
pass to the Riemannian universal covering $\widetilde{M}$ of $M$.
Indeed, by Fisher-Colbrie and Schoen result \cite{FCS} stability
is equivalent to the existence of a positive solution $u$ of
$Lu=0$ on $M$; $u$ can be lifted up by composition with the
covering projection, which is a local isometry, yielding a
positive solution of the same equation on $\widetilde{M}$.
Moreover, in this case the existence of a pole is automatically
satisfied since by \eqref{gaussequations} $M$ has non positive
sectional curvature.}
\end{Remark}
The next Lemma is a calculus exercise (see \cite{RS}).

\begin{Lemma}\label{du}
\begin{displaymath}
\mathit{If} \quad \frac{r}{\mathrm{Vol}(B_r)} \not \in
L^1(+\infty) \quad , \quad \mathit{then} \quad
\frac{1}{\mathrm{Vol}(\partial B_r)} \not \in L^1(+\infty).
\end{displaymath}
\end{Lemma}
Now we are ready to prove Theorem \ref{superfici} and Corollary \ref{corpiatto}.\\
\par
\noindent \textbf{Proof (of Theorem \ref{superfici}). } By
assumption there exists a relatively compact set $\Omega$ such
that $M\backslash \Omega$ is stable, that is,
$\lambda_1^L(M\backslash \Omega) \ge 0$. Moreover, by Gauss
equation \eqref{gaussequations} we get $2K = - |II|^2\le 0$, so
that every point of $M$ is a pole. Lemma \ref{parab} implies that
$\mathrm{Vol}(B_r) \le Cr^2$, hence
$$
\frac{r}{\mathrm{Vol}(B_r)} \not\in L^1(+\infty).
$$
From Lemma \ref{du} we obtain $(\mathrm{Vol}(\partial B_r))^{-1}
\not\in L^1(+\infty)$, and by a classical result $M$ is parabolic.
Suppose now that (\ref{supminimali}) is false, that is,
$$
\int_M |K| = \infty.
$$
Then, the function $a(x)=-2K$ satisfies all the assumptions of
Theorem \ref{spettro}, and case $(ii)$ implies that $\Delta - 2K
\equiv \Delta + |II|^2 = L$ is unstable at infinity, which is a
contradiction and concludes the proof. $\blacktriangle$\\
%
%
\par

\noindent\textbf{Proof (of Corollary \ref{corpiatto}). } If $M$ is
stable, then there exists a global positive smooth solution $u$ of
$Lu=0$. Lifting $u$ to the universal covering $\widetilde{M}$ we
deduce that $\widetilde{M}$ is a stable minimal surface with non
positive sectional curvature. By Theorem \ref{superfici},
$\widetilde{M}$ is parabolic, hence $u$ is a positive constant:
indeed,
$$
\Delta u = - |II|^2u \le 0.
$$
Equality $Lu=0$ shows that $|II|^2 \equiv 0$. Alternatively, one
can conclude as follows: by Lemma \ref{parab} we deduce $1/v
\not\in L^1(+\infty)$, where $v$ is the volume of the geodesic
spheres of $\widetilde{M}$; applying Theorem \ref{spettro}, case
$(i)$ we deduce that, if $|II|^2 \not\equiv 0$,
$\lambda_1^L(\widetilde{M})<0$, contradicting the stability
assumption. $\blacktriangle$\\
\par
\begin{Remark} \emph{Theorems \ref{superfici} and \ref{corpiatto}
can be slightly generalized to the case $\overline{\mathrm{Ricc}} \ge 0$, assuming a-priori that $M$ has a pole.
Indeed, with different techniques, by \cite{FC} there is no need to require the existence of the pole. However,
this seems to be essential in Lemma \ref{parab} in order to apply the Gauss-Bonnet theorem.\\}
\end{Remark}

\subsection{The growth of the spectral radius: proof of Theorem \ref{principale}}
We begin with a lemma. In case the volume growth is at most
exponential, by a direct application of this result we recover Do
Carmo and Zhou estimates (\ref{stimedocz1}) and
(\ref{stimedocz2}).
\begin{Lemma}\label{cor:cinque}
Suppose that
$$
\mathrm{Vol} (\partial B_r) \le f(r) \quad {\rm on} \ (R,+\infty)
$$
for some $R$ sufficiently large and some $f \in C^0
([R_0,+\infty))$. Fix $R\ge 0$.
\begin{itemize}
\item[-] If $M$ has infinite volume and
    $(\mathrm{Vol}(\partial B_r))^{-1}\not \in L^1(+\infty)$
    then
\begin{equation} \label{due:dieci}
\lambda_1^\Delta (M \backslash B_{R}) =0
\end{equation}
\item[-] If $(\mathrm{Vol}(\partial B_r))^{-1}\in
    L^1(+\infty)$, then for every $\epsilon>0$ there exists
    $T_0=T_0(\epsilon)> R$ such that
\begin{equation} \label{due:undici}
\lambda_1^\Delta (M \backslash B_{R}) \le \left\{
\inf_{t>T_0} \left[ -\frac{1}{2}
\frac{\log\int_t^{+\infty} \frac{ds}{f(s)}}{t -T_0}
\right] \right\}^2 +\epsilon
\end{equation}
\end{itemize}
\end{Lemma}
\noindent \textbf{Proof. } Set $v(r)= \mathrm{Vol}(\partial B_r)$.
We begin with the case $1/v \in L^1(+\infty)$. Let $R>0$ be
sufficiently large that
$$
R_0 > R \quad , \quad \int_{R_0}^{+\infty} \frac{ds}{v(s)} <1
$$
and let $\epsilon>0$. We define on $[R_0,+\infty)$
$$
A_\epsilon (r) =\left\{ \inf_{t>r} \left[ -\frac{1}{2}
\frac{\log\int_t^{+\infty} \frac{ds}{f(s)}}{t -r}
\right] \right\}^2+\epsilon
$$
Then, $A_\epsilon (r) \ge \epsilon$, $A_\epsilon (r)$ is
continuous and non decreasing. By Remark \ref{inL1}, $M$ has
infinite volume, thus we can apply $(v)$ of Proposition
\ref{condsuff} to obtain that (\ref{cauchy2inf}) (with
$A_\epsilon$ instead of $A$) is oscillatory. Let $z_\epsilon$ be a
locally Lipschitz solution of (\ref{cauchy2inf}), and
$R_0<T_1<T_2$ be two consecutive zeros. Define $\phi(x) =
z_\epsilon(r(x))$ on $B_{T_2} \backslash B_{T_1}$. Proceeding as
in the proof of Theorem \ref{spettro}, by the domain monotonicity
of eigenvalues we have
$$
\begin{array}{ll}
0\le \displaystyle \lambda_1^\Delta (M \backslash B_{R}) <
\lambda_1^\Delta (B_{T_2} \backslash B_{T_1}) &\le \displaystyle
\frac{\int_{B_{T_2} \backslash B_{T_1}} |\nabla
\phi|^2}{\int_{B_{T_2} \backslash B_{T_1}} \phi^2} =
\frac{\int_{T_1}^{T_2}[z_\epsilon'(r)]^2 v(r) dr}{\int_{T_1}^{T_2}z_\epsilon(r)^2
v(r)
dr} =\\[0.5cm]
\displaystyle & \displaystyle =\frac{\int_{T_1}^{T_2}A_\epsilon(r)
z_\epsilon(r)^2 v(r) dr}{\int_{T_1}^{T_2} z_\epsilon(r)^2 v(r) dr} \le
A_\epsilon(T_2) =
\\[0.5cm]
\displaystyle & \displaystyle =\left\{ \inf_{t>T_2} \left[
-\frac{1}{2} \frac{\log\int_t^{+\infty}
\frac{ds}{f(s)}}{t -T_2} \right] \right\}^2+\epsilon
\end{array}
$$
Thus we get \eqref{due:undici} with $T_0=T_2$ (note that $T_0$
depends on $\epsilon$ since $z_\epsilon(t)$ does).
\par
In case $1/v \not\in L^1(+\infty)$ and $M$ has infinite volume, by
Theorem \ref{corforte} equation (\ref{cauchy2}) is oscillatory
whenever $A(r) \ge \epsilon >0$: indeed
$$
\int_{R_0}^{+\infty}A(s)v(s)ds \ge \epsilon \int_{R_0}^{+\infty} v(s)ds
=+\infty.
$$
Thus, choosing $A_\epsilon(r)=\epsilon$ the above reasoning shows
that $\lambda_1^\Delta (M \backslash B_{R}) \le \epsilon$, and the
validity of \eqref{due:dieci} follows at once. $\blacktriangle$\\
\begin{Lemma}\label{lemdue}
In case $1/ v \in L^1(+\infty)$, the previous Lemma yields in
particular the weaker estimate
\begin{equation} \label{}
\lambda_1^\Delta (M \backslash B_{R}) \le \left\{
\liminf_{t\to +\infty} \left[ -\frac{1}{2}
\frac{\log\int_t^{+\infty} \frac{ds}{f(s)}}{t} \right]
\right\}^2 \qquad \forall \ R >0
\end{equation}
\end{Lemma}
\noindent \textbf{Proof. } This follows immediately from the next
observation: if we substitute in \eqref{due:undici} ``inf" with
the greater ``liminf", we observe that this does not depend on
$T_0(\epsilon)$. We can thus fix a particular $T_0(\epsilon)$,
compute the ``liminf"
and then let $\epsilon\to 0$. $\blacktriangle$\\
\par
\noindent \textbf{Proof (of Theorem \ref{principale}). } First, we
apply Lemma \ref{lemdue} to estimate $\lambda_1^\Delta (M
\backslash B_R)$ in case the volume growth is at most exponential.
Towards this aim suppose that $(\mathrm{Vol}(\partial B_r))^{-1}
\in L^1(+\infty)$ and that
\begin{equation} \label{due:tredice}
\mathrm{Vol} (\partial B_r) \le f(r) = \Lambda\exp\{ar^\alpha\}
\qquad 0<\alpha \le 1 \quad , \quad \Lambda,a>0
\end{equation}
Due to our choice of $\alpha$ we easily see that
$$
-\frac{1}{2}\frac{\log\int_t^{+\infty}
{\frac{ds}{f(s)}}}{t} \quad \sim \quad \frac{a}{2}
t^{\alpha -1} \qquad {\rm as \ } t \to \infty
$$
Because of this we can apply Lemma \ref{lemdue}. Hence, for every
$R\ge 0$
\begin{equation}\label{docarmo}
\lambda_1^\Delta (M \backslash B_R) \le \left\{
\begin{array}{ll}
0 \qquad & {\rm if} \ 0<\alpha < 1 \\[0.1cm]
a^2/4 \qquad & {\rm if} \ \alpha=1
\end{array}
\right.
\end{equation}
In this way we recover Do Carmo and Zhou results quoted in the
Introduction, and we also show that the estimate in Lemma
\ref{cor:cinque} is sharp. The above observations work also in
case $\mathrm{Vol}(\partial B_r) \le \Lambda \exp\{ar^\alpha
\log^\beta r\})$, with $\alpha<1,\beta \ge 0$, since it is enough
to note that
$$
\exp\{ar^\alpha \log^\beta r\} = O(\exp\{ar^{\overline{\alpha}}\})
\qquad \mathrm{for \ every} \ \overline{\alpha} > \alpha
$$
and to choose $\overline{\alpha}$ such that
$\alpha<\overline{\alpha}<1$.\par
We are left with the case $\alpha \ge 1,\beta \ge 0$. For $c>1$
and $r> R$ we define
$$
A(r) = \Big[c\Big( \frac{a\alpha}{2}\Big) r^{\alpha-1}\log^\beta
r\Big]^2.
$$
Note that $A(r)$ is monotone increasing. Moreover, Remark \ref{ft}
ensures that (\ref{cauchy2inf}) is oscillatory. Hence, proceeding
as in Lemma \ref{cor:cinque} we have for $R\ge R_0$
$$
\lambda_1^\Delta (M\backslash B_R) \le A(T_2),
$$
where $T_2(R)$ is the second zero of the solution $z$ of
(\ref{cauchy2inf}) after $R$. By Theorem \ref{teo:princ}, for
every $\varepsilon>0$ there exists $R_1(\varepsilon)$ such that,
for every $R\ge R_1$
$$
T_2(R) \le \Big[\Big(\frac{c+1}{c-1}\Big)^{\frac{2}{\alpha}}(1+
\varepsilon)\Big]R
$$
Therefore, from the monotonicity of $A(r)$ we get
$$
\lambda_1^\Delta(M\backslash B_R) \le
A\Big(\Big[\Big(\frac{c+1}{c-1}\Big)^{\frac{2}{\alpha}}(1 +
\varepsilon)\Big]R\Big) \qquad \forall R \ge R_1(\varepsilon).
$$
Inserting the value of $A(r)$, up to choosing $\varepsilon$ small
enough and $R_2 \ge R_1$ large enough we deduce that, for every
fixed $c>1$,
$$
\lambda_1^\Delta(M\backslash B_R) \le
\frac{a^2\alpha^2}{4}R^{2(\alpha-1)}\log^{2\beta} R
\Big[c^2\Big(\frac{c+1}{c-1}\Big)^{\frac{4(\alpha-1)}{\alpha}}\Big](1+2\varepsilon)
\qquad \forall \ R \ge R_2(\varepsilon).
$$
Thus, letting first $R\to +\infty$ and then $\varepsilon \to 0$,
and minimizing over all $c\in (1,+\infty)$ we finally have
\begin{equation}\label{sssss}
\limsup_{R\to +\infty}
\Bigg(\frac{\lambda_1^\Delta(M\backslash
B_R)}{R^{2(\alpha-1)}\log^{2\beta}R} \Bigg)\le
\frac{a^2\alpha^2}{4} \inf_{c\in (1,+\infty)} \Big\{c^2
\Big(\frac{c+1}{c-1}\Big)^{\frac{4(\alpha-1)}{\alpha}}\Big\}.
\end{equation}
This concludes the proof of the theorem. $\blacktriangle$\\
\par
\begin{Remark}
\emph{The infimum of the function
$$
c^2\Big(\frac{c+1}{c-1}\Big)^{\frac{4(\alpha-1)}{\alpha}}
$$
is attained by the unique positive solution $c$ of $\alpha(c+1)(c-1)= 4(\alpha-1)c$, which can be computed,
although its explicit expression is not so neat.}
\end{Remark}
\begin{Remark}
\emph{It is worth to point out that an application of
(\ref{sssss}) in case $\alpha=1$ and $\beta=0$ gives
$\lambda_1^\Delta(M\backslash B_R) \le a^2/4$, hence estimate
(\ref{sssss}) is sharp with respect to the constant appearing in
the $RHS$.}
\end{Remark}
\begin{Remark}
\emph{Proceeding as in the Introduction, one can study a model
manifold whose function $h(r)$ is of the following type:
$$
\begin{array}{ll}
h(r) = \left\{\begin{array}{ll}
\displaystyle r & r \in [0,1] \\[0.2cm]
\displaystyle \exp\left\{\frac{ar^\alpha}{m-1} \log^\beta r\right\} & r \in
[2,+\infty)
\end{array}\right. \quad ,
\end{array}
$$
for which the volume growth of geodesic spheres is
$$
\exp\{ar^\alpha \log^\beta r\}.
$$
Performing the same computations of the Introduction, one obtains
for $R$ sufficiently large
$$
\begin{array}{l}
\lambda_1^{\Delta} (M\backslash B_R) \ge K
R^{2(\alpha-1)}\log^{2\beta}R
\end{array}
$$
for some $K >0$. This shows that the estimate of Theorem
\ref{principale} is sharp even with respect to the power of the
logarithm.}
\end{Remark}
\section{Appendix}\label{appendice}
This appendix is devoted to showing existence for the Cauchy
problem \eqref{cauchy2} under general assumptions on $v(t), A(t)$.
Moreover, we prove that the zeros of such solutions, if any, are
at isolated points, and we stress a Sturm type comparison result.
In this respect, we fix $R\in (0,+\infty]$ (note that $+\infty$ is
allowed), and we assume that $v(t),A(t)$ satisfy the following set
of assumptions:
\begin{align}
&0\le A(t) \in L^\infty_{\mathrm{loc}}([0,R)) \quad , \quad A \not\equiv 0 \quad
\text{in } L^\infty_\mathrm{loc} \text{ sense} & & \tag{A1} \label{AA1}     \\
&0 \le v(t) \in L^\infty_{\mathrm{loc}}([0,R)) \quad, \quad
\frac{1}{v(t)} \in L^\infty_{\mathrm{loc}}((0,R)) \quad , \quad
\lim_{t\rightarrow 0^+} v(t)=0 & & \tag{V1} \label{VV1} \\
& \text{there exists} \ a\in (0,R) \text{ such that} \ v \text{ is
strictly increasing on} \ (0,a) & & \tag{V3} \label{VV2}
\end{align}
\begin{Proposizione}\label{esist}
Under assumptions \eqref{AA1}, \eqref{VV1}, \eqref{VV2} there
exists a locally Lipschitz function $z \in
\mathrm{Lip}_{\mathrm{loc}}([0,R))$ such that
\begin{equation}\label{quattro:uno}
\left\{
\begin{array}{l}
(v(t) z'(t))' + A(t) v(t) z(t)=0 \qquad \text{almost everywhere on
} \
(0,R)\\[0.2cm]
z(0^+)=z_0
>0
\end{array}
\right.
\end{equation}
Moreover, up to a zero-measure set $\Omega$,
$$
\lim_{\scriptsize{\begin{array}{c}
t\rightarrow 0^+\\
t \not\in \Omega \end{array} }} z'(t) = 0.
$$
If in addiction $A(t),v(t)$ are continuous on $[0,R)$, then $z \in
C^1([0,R))$ and $z'(0^+)=0$.
\end{Proposizione}
\noindent \textbf{Proof. } First, fix a sequence $T_j \uparrow R$.
We can suppose that $a \in (0,T_j)$ for every $j$, where $a$ is as
in \eqref{VV2}, and $A\not\equiv 0$ on $[0,T_j]$: the case
$A\equiv 0$ is easier and can be treated similarly. Fix
$\varepsilon \in (0,a)$, and define
$$
v_\varepsilon (t) = \left\{
\begin{array}{ll}
v(\varepsilon) &{\rm on} \ (0,\varepsilon]\\[0.2cm]
v(t) &{\rm on} \ [\varepsilon,R)
\end{array}
\right.
$$
Then,
\begin{equation}\label{quattro:due}
k_\varepsilon (t,s)=-A(s)v_\varepsilon (s) \int_s^t
\frac{dx}{v_\varepsilon (x)}
\end{equation}
belongs to $L^\infty_\mathrm{loc}([0,R) \times  [0,R))$.  Thus, by
standard theory (one can consult chapter IX of \cite{KF}),
Volterra integral equation of the second type
\begin{equation}\label{quattro:tre}
w(t)=z_0+\int_0^t k_\varepsilon (t,s) w(s) ds,
\end{equation}
restricted to every interval $[0,T_j]$ (where the kernel
$k_\varepsilon(t,s)$ is bounded), admits a unique solution
$z_{\varepsilon,j} \in L^2((0,T_j))$. From \eqref{quattro:due},
using integration by parts applied to the integrable function
$-A(s)v_\varepsilon(s)z_{\varepsilon,j}(s)$ and to the absolutely
continuous one
$$
\int_s^t \frac{dx}{v_\varepsilon(x)}
$$
We see that $z_{\varepsilon,j}$ satisfies
\begin{equation} \label{quattro:quattro}
z_{\varepsilon,j}(t)= z_0 -\int_0^t \frac{1}{v_\varepsilon (s)}
\left\{ \int_0^s A(x) v_\varepsilon (x) z_{\varepsilon,j}(x) dx
\right\} ds
\end{equation}
on $[0,T_j]$. This shows that $z_{\varepsilon,j}(t)$, being an
integral function, is absolutely continuous on $[0,T_j]$ (hence,
almost everywhere differentiable), and its derivative is almost
everywhere
$$
\frac{1}{v_\varepsilon (t)}  \int_0^t A(x) v_\varepsilon (x)
z_{\varepsilon,j}(x) dx \in L^\infty([0,T_j]).
$$
Therefore, $z_{\varepsilon,j}(t)$ is a Lipschitz function on
$[0,T_j]$. By the uniqueness of solutions of \eqref{quattro:tre},
we deduce that, when $j<k$, $z_{\varepsilon,k}$ restricted to
$[0,T_j]$ coincides with $z_{\varepsilon,j}$. Hence, we can
construct a locally Lipschitz solution $z_{\varepsilon}(t)$
defined on the whole $[0,R)$. What we want to prove is that, for
every $T_j$, the family $\{ z_\varepsilon \}_{\varepsilon \in
(0,a)}$ is equibounded and equi-Lipschitz in $C^0 ([0,T_j])$. For
the ease of notation, from now on we omit the subscript $j$ and we
consider the problem on $[0,T] \subset [0,R)$. We observe that,
because of \eqref{VV2} and \eqref{AA1}, for $0 \le s \le t \le a$
we have
\begin{equation}\label{quattro:cinque}
|k_\varepsilon (t,s)| \le \|A\|(t-s) \le \|A\|a
\end{equation}
where $\|A\|= \|A\|_{L^\infty([0,T])}$. Next, we consider the case
$0 \le s \le a < t \le T$. Because of \eqref{VV1}, on $[a,T]$
$v^{-1}$ is bounded. We indicate with $\|v^{-1}\|$ the
$L^\infty$-norm of $v^{-1}(t)$ on $[a,T]$, and with $\|v\|$ the
$L^\infty$-norm of $v(t)$ on the whole $[0,T]$. It follows that
$$
\begin{array}{ll}
\displaystyle |k_\varepsilon (t,s)| &= \displaystyle A(s)
v_\varepsilon (s)\left\{ \int_s^a \frac{dx}{v_\varepsilon
(x)}+\int_a^t \frac{dx}{v_\varepsilon(x)}   \right\} \le
\|A\| \Big(a + v_\varepsilon (s)\|v^{-1}\|T \Big)
\\[0.4cm]
 & \displaystyle \le \|A\| \left( a + \|v\|\|v^{-1}\|T \right)
\end{array}
$$
It remains to consider the case $0<a \le s\le t \le T$. In this
case we obtain
$$
|k_\varepsilon (t,s)| \le A(s) v_\varepsilon (s) \|v^{-1}\|T \le \|A\|\|v\|\|v^{-1}\|T.
$$
Therefore, there exists $L=L(T,a) >0$ such that
\begin{equation}\label {quattro:sette}
\sup_{\varepsilon \in (0,a)} \left( \sup_{0 \le s \le t \le T}
|k_\varepsilon (t,s)|  \right) \le L
\end{equation}
Using \eqref{quattro:sette} into \eqref{quattro:tre} we have
$$
|z_\varepsilon (t)| \le z_0 +L \int_0^t |z_\varepsilon
(s)| ds \qquad \forall \ t \in [0, T]
$$
So that, applying Gronwall lemma on the continuous function $|z_\varepsilon (t)|$, we conclude
\begin{equation}\label{quattro:otto}
|z_\varepsilon (t)| \le z_0 {\rm e}^{Lt} \le z_0 {\rm e}^{LT}
\qquad {\rm on} \ [0,T]
\end{equation}
This shows equiboundedness of the family $\{ z_\varepsilon
\}_{\varepsilon \in (0,a)}$. To show equicontinuity we
differentiate \eqref{quattro:quattro} to obtain
\begin{equation}\label{quattro:nove}
z_\varepsilon ' (t) = -\frac{1}{v_\varepsilon (t)} \int_0^t
A(x) v_\varepsilon (x) z_\varepsilon (x) dx \qquad \text{almost everywhere on }
[0,T].
\end{equation}
We set
$$
H(\varepsilon,t)=\frac{1}{v_\varepsilon (t)} \max_{s \in [0,t]}
A(s) v_\varepsilon (s)
$$
If $0 \le t \le a$, because of \eqref{AA1} and \eqref{VV2} we have
$$
H(\varepsilon, t)\le \|A\|
$$
If $a < t \le T$, since $\varepsilon \in (0,a)$, $v_\varepsilon
(t)=v(t)$ and therefore
$$
H(\varepsilon,t) \le \|A\|
\frac{\|v_\varepsilon\|_{L^\infty([0,t])}}{v(t)} \le \|A\|
\|v^{-1}\| \|v_\varepsilon\|_{L^\infty([0,t])} \le
\|A\|\|v^{-1}\|\|v\|,
$$
where the last inequality is an immediate consequence of
\eqref{VV2} and the definition of $v_\varepsilon(t)$. Summarizing,
there exists $M=M(T,a)>0$ such that
$$
\sup_{\varepsilon \in (0,a)} H (\varepsilon,t) \le M \qquad
\text{a.e. on } [0,T].
$$
From \eqref{quattro:nove} it follows that
$$
|z_\varepsilon ' (t)| \le M \int_0^t |z_\varepsilon (x)| dx \qquad
\text{a.e. on } [0,T]
$$
and thus, from \eqref{quattro:otto},
\begin{equation}\label{stimaderivata}
|z_\varepsilon ' (t)| \le z_0 MT {\rm e}^{LT} \qquad \text{a.e. on
} [0,T]
\end{equation}
This shows that $\{ z_\varepsilon \}_{\varepsilon \in (0,a)}$ is
equi-Lipschitz on every compact subset $[0,T]\subset [0,R)$. By
the Ascoli-Arzel\'a theorem, the set $\{ z_\varepsilon
\}_{\varepsilon \in (0,a)}$ is relatively compact in $C^0
([0,T])$. Therefore, there exists a sequence $\varepsilon_n \to 0$
such that $z_{\varepsilon_n}$ converges uniformly to a Lipschitz
function $z$ on $[0,T]$. A Cantor diagonal argument on the
exhaustion $[0,T_j] \uparrow [0,R)$ yields a sequence
$z_{\varepsilon_n}$ which converges locally uniformly to a locally
Lipschitz function $z$ on
$[0,R)$.\\
Clearly, $v_{\varepsilon_n} \to v$ in $L^\infty([0,R))$. If we set
$$
r_\varepsilon (t) =\frac{1}{v_\varepsilon (t)} \int_0^t {A(s)
v_\varepsilon (s) z_\varepsilon (s) ds}
$$
using \eqref{quattro:nove} and \eqref{stimaderivata} we see that
$r_{\varepsilon_n}$ is locally a bounded sequence of
$L^\infty_\mathrm{loc}$-functions converging pointwise to
$$
r(t) =\frac{1}{v (t)} \int_0^t {A(s) v (s) z(s) ds} \qquad
\text{a.e. on }  [0,R)
$$
By the dominated convergence theorem $r_{\varepsilon_n} \to r$ in
$L^1 ((0,t])$ $\forall \ t \in (0,R)$. Hence, for every $t \in
[0,R)$,
$$
\lim_{n \to +\infty} \int_0^t \frac{ds}{v_{\varepsilon_n} (s)}
\left\{ \int_0^s A(x) v_{\varepsilon_n} (x)
z_{\varepsilon_n}(x) dx \right\} = \int_0^t   \frac{ds}{v
(s)} \left\{ \int_0^s A(x) v(x) z (x) dx  \right\}
$$
Because of \eqref{quattro:quattro}  it follows  that $z$ satisfies
the integral equation
\begin{equation}\label{sol_int}
z(t) = z_0 - \int_0^t \frac{1}{v(s)} \left\{ \int_0^s A(x) v(x)
z (x) dx \right\} ds,
\end{equation}
hence the Cauchy problem \eqref{quattro:uno}. Note that, in case
$v(t),A(t)$ are also continuous, from \eqref{sol_int} we deduce
that $z(t) \in C^1((0,R))$. Because of \eqref{VV2}, for $t \in
(0,a]$ we have
$$
|z'(t)| = \left|\frac{1}{v(t)} \left\{ \int_0^t A(s) v(s) z (s) ds
\right\} \right|
 \le \int_0^t {A(s) |z(s)|} ds \qquad \text{ almost
everywhere}
$$
so that, up to a zero-measure set $\Omega$, $z'(t) \to 0$ as $t
\to 0^+$, $t\not\in \Omega$. In case $v(t),A(t)$ are continuous,
the above inequality is everywhere valid and shows that $z(t) \in
C^1([0,R))$ with $z'(0^+)=0$. This concludes the proof.
$\blacktriangle$\\
\begin{Remark} \emph{With the same technique (but a simpler proof) we
can provide existence of a locally Lipschitz solution of problem
\eqref{cauchy2inf} when \eqref{AA1}, \eqref{VV1} are met on
$[t_0,R)$, for some $t_0>0$. Note that $1/v$ is required to be
bounded also in a neighborhood of $t_0$.}\\
\end{Remark}

\begin{Proposizione}\label{isolati}
Assume \eqref{AA1} and \eqref{VV1}. Then, the zeros of every
Locally Lipschitz solution $z(t)$ of \eqref{quattro:uno}, if any,
are at isolated points of $[0,R)$.
\end{Proposizione}

\noindent \textbf{Proof. } Let
$$
y(t) = -\frac{v(t)z'(t)}{z(t)}.
$$
Since $z \in \mathrm{Lip}_\mathrm{loc}([0,R))$, $y(t)$ is at least
locally Lipschitz on compact sets of $[0,R)\backslash \{t :
z(t)=0\}$. This follows since $(vz')' = -Avz \in
L^\infty_\mathrm{loc}([0,R))$, hence $vz'$ is locally Lipschitz.
Differentiating and using \eqref{quattro:uno} we get
$$
y'(t) = A(t)v(t) + \frac{y^2(t)}{v(t)} \qquad \text{almost everywhere,}
$$
hence $y(t)$ is increasing on its domain. Assume that $t_0\in
(0,R)$ is a zero of $z(t)$ (note that $z_0>0$). First, we prove
that
\begin{equation}\label{limitiy}
\lim_{t\rightarrow t_0^\pm} y(t) = \mp \infty
\end{equation}
Indeed, both limits exist by monotonicity. Indicating with $L^\pm$
the two limits, if by contradiction $L^- <+\infty$ (analogously
for $L^+>-\infty$) then necessarily
$$
v(t_0)z'(t_0)= \lim_{t\rightarrow t_0}v(t)z'(t) = \lim_{t\rightarrow
t_0^-} v(t)z'(t) = -z(t_0)L^-=0.
$$
Therefore, $z(t)$ should solve
\begin{equation}
\left\{
\begin{array}{l}
(v(t) z')' + A(t) v(t) z=0 \qquad \text{almost everywhere on } \
(0,R)\\[0.2cm]
z(t_0)=0 \qquad , \qquad v(t_0)z'(t_0) = 0.
\end{array}
\right.
\end{equation}
In other words, $z(t)$ should be a locally Lipschitz solution of
Volterra integral problem
\begin{equation}\label{integvol}
z(t)= - \int_{t_0}^t \frac{1}{v(s)} \left\{ \int_{t_0}^s A(x) v(x)
z (x) dx \right\} ds = -\int_{t_0}^t \left[ A(s)v(s)\int_s^t \frac{dx}{v(x)}\right] z(s)ds,
\end{equation}
where the last inequality follows integrating by parts. Since
$v(t)$ is bounded away from zero on compact sets of $(0,R)$, the
kernel of Volterra operator is locally bounded. Therefore,
\eqref{integvol} has a unique local solution, which is necessarily
$z\equiv 0$ on every $[T_1,T_2]\subset (0,R)$. This contradicts
$z(0^+)= z_0>0$ and proves \eqref{limitiy}. Now, if there exists
$\{t_k\}$ such that $z(t_k)=0$ and $t_k \rightarrow t_0$, every
neighborhood of $t_0$ should contain points $t_k$ such that
$\lim_{t\rightarrow t_k^\pm} y(t) = \mp \infty$, and this clearly
contradicts the fact that $y(t)$ has both left and right limits in
$t_0$. $\blacktriangle$\\

\begin{Proposizione}\label{comparison}
Assume \eqref{VV1}, and let $A_1, A_2$ satisfy \eqref{AA1} and
$A_1 \ge A_2$ a.e. on $[0,R)$. Suppose that $z_i(t)$, $i\in
\{1,2\}$, is a locally Lipschitz solution of \eqref{quattro:uno}
with $A(t)=A_i(t)$. Fix $T \le R$ such that $z_1(t),z_2(t)> 0$ on
$[0,T)$. Then $z_2(t) \ge z_1(t)$ on $[0,T)$.
\end{Proposizione}
\noindent \textbf{Proof. } We consider the locally Lipschitz
function $F= (vz_1')z_2 - (v z_2')z_1$. Differentiating we obtain
$$
F' = \Big(vz_1'z_2 - v z_2'z_1\Big)' = z_2(-A_1vz_1) -
z_1(-A_2vz_2) = (A_2-A_1)vz_1z_2 \le 0
$$
almost everywhere on $[0,T)$. This shows that $F$ is non
increasing. From $F(0^+)=0$ we argue $F\le 0$ and therefore
$vz_1'z_2 \le vz_2'z_1$. By \eqref{VV1} we deduce that $v$ is
essentially bounded from below with a positive constant on compact
sets of $(0,T)$, thus $z_1'z_2 \le z_2'z_1$ almost everywhere.
Hence
$$
\left(\frac{z_1}{z_2}\right)' = \frac{z_1'z_2-z_2'z_1}{z_2^2} \le
0 \qquad \text{almost everywhere on } (0,T).
$$
Since $z_1(0)/z_2(0) = 1$ we conclude $z_1(t)\le z_2(t)$ on
$[0,T)$.

\vspace{0.5cm}

\noindent \textbf{Acknowledgements. } The authors express their
gratitude to the referee for a very careful reading of the
manuscript and for the many useful observations which led to
substantial improvements.

\end{document}